\documentclass[reqno,11pt]{amsart}
\usepackage{amsmath,amssymb,amsthm, verbatim}
\usepackage{url}
\usepackage[usenames, dvipsnames]{color}

\usepackage[letterpaper,hmargin=1in,vmargin=1in]{geometry}

\usepackage{graphicx}
\usepackage{enumerate,pinlabel}
\usepackage{mathrsfs,graphicx,url}
\usepackage[usenames,dvipsnames]{xcolor}
\usepackage[colorlinks=true,linkcolor=Black,citecolor=Black, urlcolor=Black, pagebackref]{hyperref}

\numberwithin{equation}{section}

\theoremstyle{plain}
\newtheorem{theorem}{Theorem}
\newtheorem{conjecture}{Conjecture}
\newtheorem*{theorem*}{Theorem}
\newtheorem*{lemma*}{Lemma}
\newtheorem{lemma}[theorem]{Lemma}
\theoremstyle{remark}
\newtheorem{remark}[theorem]{Remark}

\theoremstyle{definition}

\theoremstyle{remark}

\newcommand{\NN}{\mathbb{N}}
\newcommand{\R}{\mathbb{R}}

\newcommand\supp{\mathop{\rm supp}}

\newcommand\Tr{\mathop{\rm Tr}}

\newcommand*{\defeq}{\mathrel{\vcenter{\baselineskip0.5ex \lineskiplimit0pt

                     \hbox{\scriptsize.}\hbox{\scriptsize.}}}%
                     =}

\title{Asymptotics for Riesz means of integrable polygons}

\begin{document}

\author{Amanda Criner}
\author{Gustav M\r{a}rdby}
\author{Jacob Shapiro}

\address{Air Force Research Laboratory, Wright-Patterson AFB, OH, USA}

\address{Mathematical Sciences, Chalmers University of Technology and University of Gothenburg, 412 96 Gothenburg, Sweden}
\email{mardby@chalmers.se}

\address {Department of Mathematics, University of Dayton, Dayton, OH 45469-2316, USA}
\email{jshapiro1@udayton.edu}

\begin{abstract}
We prove complete asymptotic expansions for Riesz means of order $\gamma > 0$ for the Dirichlet and Neumann Laplacians on integrable polygons. The coefficients of the first three terms are the familiar ones involving the area, perimeter, and corner angles of the polygon. The coefficients of the lower-order terms are given by convergent trigonometric series. %The proof is based on the heat trace identities for integrable polygons obtained by Rowlett and the third author \cite{maro25}, together with large-argument asymptotics for Bessel functions. 
The proof is based on heat trace identities for integrable polygons, together with large-argument asymptotics for Bessel functions.
\end{abstract}
\maketitle 
%\author

\section{Introduction}

In this paper we investigate the classical problem of the asymptotic distribution of Laplace eigenvalues in the setting of integrable polygons in the plane. By an integrable polygon we mean a polygon whose interior angles are all of the form $\pi/n$, where $n \geq 2$ is an integer. In particular, every interior angle is at most $\pi/2$, so every integrable polygon is convex.

The class of integrable polygons is small. Indeed, since the sum of the interior angles of a convex polygon with $N$ sides is $(N-2)\pi$, while each interior angle of an integrable polygon is at most $\pi/2$, we have $(N-2)\pi \leq N\pi/2$, or $N \leq 4$. It follows that the only integrable polygons are rectangles, equilateral triangles, isosceles right triangles, and hemi-equilateral (30-60-90) triangles. %the polygonal domains 
These polygons are precisely the ones which strictly tessellate the plane  \cite{mc08}. %are precisely rectangles, equilateral triangles, isosceles right triangles, and hemi-equilateral triangles. 
%Therefore, a polygon is integrable if and only if it strictly tessellates the plane.

The study of eigenvalue asymptotics goes back to the celebrated work of Weyl \cite{we11}. If $\Omega \subseteq \mathbb R^2$ is a polygon and $\{\lambda_j\}$ denotes the eigenvalues of the Laplacian $-\Delta$ on $\Omega$, repeated according to multiplicities and subject to either the Dirichlet or Neumann boundary condition, the eigenvalue counting function is
\begin{equation*}
    N(E) \defeq \#\{j : \lambda_j \leq E\}.
\end{equation*}
In this setting, Weyl's law says
\begin{equation} \label{weyl}
    N(E) = \frac{|\Omega|}{4\pi} E + o(E), \quad E \to \infty,
\end{equation}
where $|\Omega|$ denotes the Lebesgue measure or area of the domain. For the Dirichlet boundary condition, \eqref{weyl} holds under the assumption that $\Omega \subseteq \R^2$ is open with finite measure \cite{ro72}. For the Neumann Laplacian, a sufficient condition is that $\Omega \subseteq \mathbb R^2$ has the extension property \cite{flw23}. This means that there exists a bounded linear operator $\mathcal E : H^1(\Omega) \to H^1(\mathbb R^2)$ such that $\mathcal E u = u$ almost everywhere on $\Omega$ for every $u \in H^1(\Omega)$. Here $H^1$ denotes the usual first order Sobolev space.

It is natural to ask whether the remainder in Weyl's law can be improved. For domains with smooth boundary, Seeley \cite{se78} showed that the counting function satisfies the sharper estimate
\begin{equation*}
    N(E)=\frac{|\Omega|}{4\pi}E+O(E^{1/2}), \quad E\to \infty,
\end{equation*}
for both the Dirichlet and Neumann boundary conditions. Ivrii \cite{iv80} proved, under an additional dynamical condition on the billiard flow, the two-term asymptotic
\begin{equation*}
    N(E)=\frac{|\Omega|}{4\pi}E \pm \frac{|\partial \Omega|}{4\pi}E^{1/2}+o(E^{1/2}), \quad E\to \infty.
\end{equation*}
Here, $|\partial \Omega|$ is the perimeter of $\Omega$; the minus sign corresponds to the Dirichlet boundary condition and the plus sign to the Neumann boundary condition. These results hold more generally in dimension $d \ge 2$, with the second term involving the $(d-1)$-dimensional surface measure of the boundary. In our polygonal setting one expects lower order contributions coming from the corners.

The study of refined Weyl asymptotics for domains with corners goes back at least to the work of Fedosov. In the plane \cite{fe63}, Fedosov showed that, for a polygonal domain $\Omega$
with either the Dirichlet or Neumann boundary condition, the eigenvalue counting function in frequency, $n(\nu) \defeq \#\{j : \lambda_j \le \nu^2\}$, has the expansion
\begin{equation*}
    n(\nu) = \frac{|\Omega|}{4\pi}\nu^{2}+O(\nu), \quad \nu \to \infty.
\end{equation*}
He also showed that its integral satisfies
\begin{equation} \label{fedosov integral}
    \int_{0}^{\nu}n(t)\,dt = \frac{|\Omega|}{12\pi}\nu^{3} \pm \frac{|\partial \Omega|}{8\pi}\nu^{2} + O(\nu),
\end{equation}
where the minus sign corresponds to the Dirichlet boundary condition and the plus sign
to the Neumann boundary condition. His strategy was based on applying a Fourier Tauberian theorem of Mar\v{c}enko \cite{ma55} to the short-time expansion of the wave trace. Fedosov generalized this result to a $d$-term expansion for polyhedra in dimension $d \ge 2$ with remainder $O(\nu^{d-1})$ \cite{fe64}.

Rather than the counting function itself, in this paper we consider its integral or averaged version (in the spirit of \eqref{fedosov integral}), the so-called \textit{Riesz mean}. It is usually defined as the trace of the operator $(-\Delta - E)_-^\gamma$, where $E,\gamma > 0$ and $x_- \defeq \max(0,-x)$. In particular, we study
\begin{equation*}
    \Tr(-\Delta - E)_-^\gamma = \sum_{\lambda_j \le E} (E - \lambda_j)^\gamma %= \int_0^E N(t)dt
    , \qquad E \to \infty.
\end{equation*}
%Our main result is a four-term asymptotic expansion for the Riesz mean for the Dirichlet and Neumann Laplacians on integrable polygons. In addition to area, perimeter, and corner contributions, \red{the expansion contains explicit almost periodic terms of successively lower orders.}
Our main result is a complete asymptotic expansion for the Riesz mean associated with the Dirichlet and Neumann Laplacians on integrable polygons. In addition to the area, perimeter, and corner contributions, the expansion contains trigonometric series of successively lower orders. 

\begin{theorem} \label{Riesz mean thm}
    Let $\gamma > 0$ and let $\Omega$ be an integrable polygon with interior angles $\theta_1,\dots,\theta_N$. Let $0 \leq \lambda_0 < \lambda_1 \leq \lambda_2 \leq \cdots$ be the eigenvalues of the nonnegative Laplacian $-\Delta$ on $\Omega$ with either the Dirichlet or Neumann boundary condition, repeated according to their multiplicities. %Then, the Riesz mean $R_\Omega^\gamma(E) \defeq \sum_{\lambda_j \leq E} (E - \lambda_j)^\gamma$ admits the asymptotic expansion
    %\begin{equation} \label{Riesz mean asymptotic}
        %\begin{split}
            %R_\Omega^\gamma(E) &\sim \frac{|\Omega|}{4\pi(\gamma + 1)}E^{\gamma + 1} \pm \frac{|\partial\Omega|}{8\sqrt{\pi}}\frac{\Gamma(\gamma + 1)}{\Gamma(\gamma + 3/2)}E^{\gamma + 1/2} + E^\gamma\sum_{i=1}^N \frac{\pi^2 - \theta_i^2}{24\pi\theta_i} \\
            %&+ \sum_{k = 1}^\infty F_{\Omega,k}^\gamma(\sqrt{E})E^{(\gamma+1)/2-k/4}, \quad E \to \infty.
        %\end{split}
    %\end{equation}
    Then, for every $K \geq 1$, the Riesz mean $R_\Omega^\gamma(E) \defeq \sum_{\lambda_j \leq E} (E - \lambda_j)^\gamma$ satisfies
    \begin{equation} \label{Riesz mean asymptotic}
        \begin{split}
            R_\Omega^\gamma(E) &= \frac{|\Omega|}{4\pi(\gamma+1)}E^{\gamma+1} \pm \frac{|\partial\Omega|}{8\sqrt{\pi}} \frac{\Gamma(\gamma+1)}{\Gamma(\gamma+3/2)} E^{\gamma+1/2} + E^\gamma \sum_{i=1}^N
            \frac{\pi^2-\theta_i^2}{24\pi\theta_i} \\
            &+ \sum_{k=1}^K F_{\Omega,k}^\gamma(\sqrt E) E^{(\gamma+1)/2-k/4} + O(E^{(\gamma+1)/2 - (K+1)/4}), \quad E \to \infty.
        \end{split}
    \end{equation}
    Here, the sign of the second term is negative for the Dirichlet boundary condition and positive for the Neumann boundary condition. Moreover, for each $k \geq 1$, the coefficient function $F_{\Omega,k}^\gamma$ is given by convergent trigonometric series; see \eqref{eq: Riesz mean rectangle asymptotic} and \eqref{eq: Riesz mean equilateral asymptotic}-\eqref{eq: Riesz mean isosceles asymptotic}. If $\gamma > 1/2$, all of these series converge absolutely and the coefficient functions are bounded. If $0 < \gamma \leq 1/2$, all of these series except the double series defining $F_{\Omega,1}^\gamma$ converge absolutely; the latter converge conditionally under rectangular summation, and $F_{\Omega,1}^\gamma$ is unbounded.
\end{theorem}

Theorem \ref{Riesz mean thm} should be compared with the recent results of Frank and Larson for Riesz means on general polygons \cite[Theorems 6.2 and 6.3]{frla25}. Their results apply to $\Tr(-\Delta-E)^\gamma_-$ for $\gamma \geq 0$. In the Dirichlet case, they obtain a non-asymptotic expansion, valid for all $E>0$, whose terms are the area, perimeter, and corner contributions, with an explicit bound on the remainder. In the Neumann case, they obtain the analogous asymptotic statement as $E \to \infty$ with remainder $O(E^{(\gamma+1)/2})$. Their approach uses a Tauberian theorem for the Laplace transform to convert the heat trace expansions for polygons of van den Berg and %S. 
Srisatkunrajah \cite{vansri88} into Riesz mean asymptotics and remainder estimates. 
For the integrable polygons, %and $\gamma > 1/2$, 
Theorem \ref{Riesz mean thm} strengthens these results by replacing the $O(E^{(\gamma+1)/2})$ remainder with a complete asymptotic expansion whose %leading term is of order $E^{\gamma/2+1/4}$.
first term is $F_{\Omega,1}^\gamma(\sqrt E)E^{\gamma/2+1/4}$.

Along with \cite{frla25}, the present paper is motivated by the recent work of the second author and Rowlett \cite{maro25}, who studied other spectral invariants of integrable polygons, including formulas for spectral zeta functions, zeta-regularized determinants, and sharp short-time heat trace asymptotics. The eigenvalues of integrable polygons are known explicitly; this tractability underlies the heat trace formulas obtained in \cite{maro25}. Our starting point is these formulas evaluated at complex time. We integrate the identities against the inverse Fourier transform of $\lambda \mapsto (|\lambda|-E)_-^\gamma$, to obtain a formula for the Riesz mean in the limit as the real part of the time tends to zero.
%To proceed, we interchange a limit with an infinite sum; the condition $\gamma > 1/2$ justifies this exchange. We write the resulting terms using Bessel functions and then apply their large-argument asymptotic expansion as $E \to \infty$. The threshold $\gamma > 1/2$ is worth noting because it is different from the division in \cite{frla25} between the ranges $\gamma \ge 1$ and $0 < \gamma < 1$, which are treated there by separate arguments. Here, the condition $\gamma > 1/2$ arises in the treatment of the double sums appearing in the heat trace formulas. After converting the sums into Bessel series, the single sums are absolutely convergent for every $\gamma > 0$, whereas the double sums are absolutely convergent precisely when $\gamma > 1/2$. Our proof relies on this absolute convergence to justify the interchange of the limit $\sigma \to 0^+$ with the infinite sums. It is natural to ask whether these formulas continue to hold for $0 < \gamma \leq 1/2$, with the double sums interpreted, for instance, by radial summation. Interestingly, the threshold $\gamma = 1/2$ coincides with the critical Bochner-Riesz index $(d-1)/2$ in dimension $d = 2$ \cite{ks19}. 
We compute the limits using Bessel functions and then apply their large-argument asymptotic expansion as $E \to \infty$. The threshold $\gamma=1/2$ plays an important role in the convergence of the resulting series. It arises from the double sums appearing in the heat trace formulas. After converting the sums into Bessel series, all resulting series are absolutely convergent for every $\gamma > 0$, except for the leading double series, which are absolutely convergent if and only if $\gamma > 1/2$. For $0 < \gamma \leq 1/2$, these double series are no longer absolutely convergent, but we show that they converge conditionally under rectangular summation. This allows us to extend the Riesz mean expansion to the full range $\gamma > 0$. The threshold $\gamma = 1/2$ is worth noting because it differs from the division in \cite{frla25} between the ranges $\gamma \geq 1$ and $0 < \gamma < 1$, which are treated there by separate arguments. Interestingly, the threshold $\gamma = 1/2$ coincides with the critical Bochner--Riesz index $(d-1)/2$ in dimension $d = 2$ \cite{ks19}.

We conjecture that the asymptotic expansion in Theorem \ref{Riesz mean thm} extends to arbitrary polygons. There are several reasons to believe this. First, as mentioned above, Frank and Larson \cite[Theorems 6.2 and 6.3]{frla25} proved that the first three terms in \eqref{Riesz mean asymptotic} hold for all polygons. Second, by van den Berg and Srisatkunarajah \cite{vansri88}, the Dirichlet heat trace of a polygon has a three-term asymptotic expansion followed by an exponentially small remainder. The analogous statement is also known for the Neumann heat trace \cite{ma26}. Motivated by these results, we make the following conjecture.
\begin{conjecture}
    Let %$\gamma > 1/2$ 
    $\gamma > 0$
    and let $\Omega$ be a polygon with interior angles $\theta_1,\dots,\theta_N$. Let $R_\Omega^\gamma(E)$ denote the Riesz mean associated with either the Dirichlet or Neumann Laplacian on $\Omega$. Then there exist %bounded 
    functions
    \begin{equation*}
        F_{\Omega,k}^\gamma:[0,\infty)\to \R, \quad k \geq 1,
    \end{equation*}
    %such that
    %\begin{align*}
        %R_\Omega^\gamma(E) &\sim \frac{|\Omega|}{4\pi(\gamma+1)}E^{\gamma+1} \pm \frac{|\partial\Omega|}{8\sqrt{\pi}} \frac{\Gamma(\gamma+1)}{\Gamma(\gamma+3/2)} E^{\gamma+1/2} + E^\gamma\sum_{i=1}^N
        %\frac{\pi^2-\theta_i^2}{24\pi\theta_i} \\
        %&+ \sum_{k=1}^{\infty} F_{\Omega,k}^\gamma(\sqrt E) E^{(\gamma+1)/2-k/4}, \quad E \to \infty.
    %\end{align*}
    such that, for every $K\geq1$,
    \begin{align*}
        R_\Omega^\gamma(E) &= \frac{|\Omega|}{4\pi(\gamma+1)}E^{\gamma+1} \pm \frac{|\partial\Omega|}{8\sqrt{\pi}} \frac{\Gamma(\gamma+1)}{\Gamma(\gamma+3/2)} E^{\gamma+1/2} + E^\gamma\sum_{i=1}^N \frac{\pi^2-\theta_i^2}{24\pi\theta_i} \\
        &+ \sum_{k=1}^K F_{\Omega,k}^\gamma(\sqrt{E}) E^{(\gamma+1)/2-k/4} + O\big(E^{(\gamma+1)/2-(K+1)/4}\big), \quad E \to \infty.
    \end{align*}
\end{conjecture}

For the integrable polygons, Theorem \ref{Riesz mean thm} proves this and shows, in addition, that the coefficient functions are given by convergent trigonometric series.

\medskip

\noindent \textbf{Acknowledgments:} We are grateful to Simon Larson for discussions, guidance and assistance throughout this work. We would also like to thank Julie Rowlett and Kiril Datchev for initial discussions. A.K. Criner and J.Z. Shapiro's work was partially supported by the Air Force Research Laboratory and the Air Force Office of Scientific Research, under Lab Task 21RXCOR069. A.K. Criner and J.Z. Shapiro’s collaboration would not have been possible without the support of Phil Buskohl. 

\medskip

\noindent \textbf{AI usage disclosure:} The authors made use of ChatGPT for mathematical discussions and editorial assistance, particularly in the development of the arguments in Appendix \ref{double sum bounds appendix}. All mathematical statements and proofs were independently checked, revised, and written in their final form by the authors.

\section{Riesz mean asymptotics for rectangles} \label{prelim section}
The eigenvalues of $-\Delta$ on a rectangle of side lengths $a$ and $b$ with the Dirichlet boundary condition are 
\begin{equation*}
\lambda^D_{m,n}(\square)= \pi^2\bigg(\frac{m^2}{a^2} + \frac{n^2}{b^2} \bigg), \quad m, n \geq 1.
\end{equation*}
The multiplicity of an eigenvalue $\lambda$ is the number of pairs $(m,n)$ such that $\lambda = \lambda^D_{m,n}(\square)$. For the Neumann boundary condition, the eigenvalues $\lambda^N_{m,n}(\square)$ and their multiplicities are the same, except now $m$ and $n$ range over the nonnegative integers. 

We begin with an identity for the heat trace %$H_\square(t)$ 
of a rectangle, obtained in \cite[Section 2]{maro25}:
\begin{equation*}
    \begin{split}
        H_\square(t) := \sum_{j = 0}^\infty e^{-t\lambda_j} &= \frac{ab}{4\pi t} \pm \frac{a +b}{4 \sqrt{\pi t}} + \frac{1}{4} \pm \frac{a}{2 \sqrt{\pi t}} \sum_{m =1}^{\infty} e^{-a^2m^2/t} \pm \frac{b}{2 \sqrt{\pi t}} \sum_{n=1}^\infty e^{-b^2n^2/t} \\
        &+ \frac{ab}{2\pi t} \sum_{m =1}^{\infty} e^{-a^2m^2/t}  + \frac{ab}{2 \pi t} \sum_{n=1}^\infty e^{-b^2n^2/t} + \frac{ab}{\pi t} \sum_{m,n \geq 1} e^{-(a^2m^2 +b^2n^2)/t}, \quad t > 0.
        \end{split}
\end{equation*}
For every term where we have $\pm$, we use plus for the Neumann boundary condition and minus for the Dirichlet boundary condition. Both sides analytically continue to a complex argument $z$ with $\operatorname{Re}(z)>0$, where the branch of $\sqrt{z}$ is chosen so that $\operatorname{Re}(\sqrt{z})>0$. Thus,
\begin{equation} \label{complex heat trace id}
    \sum_{j = 0}^\infty e^{-(\sigma + it) \lambda_j}  = \frac{ab}{4\pi( \sigma + it)} \pm \frac{a + b}{4 \sqrt{\pi(\sigma + it)}} + \frac{1}{4} + \mathcal{R}_\square(\sigma + it), \quad \sigma > 0, \quad t\in \R.
\end{equation}
Here, $\mathcal{R}_\square(\sigma + it)$ is the remainder that contains any term involving an infinite sum. In particular,
\begin{equation} \label{R of sigma plus it}
    \begin{split}
        \mathcal{R}_\square(\sigma + it) &= \frac{ab}{2\pi(\sigma + it)} \sum_{m=1}^\infty e^{-a^2m^2/(\sigma + it)} \pm \frac{a}{2\sqrt{\pi(\sigma + it)}}\sum_{m=1}^\infty e^{-a^2m^2/(\sigma + it)} \\
        &+ \frac{ab}{2\pi(\sigma + it)} \sum_{n=1}^\infty e^{-b^2n^2/(\sigma + it)} \pm \frac{b}{2\sqrt{\pi(\sigma + it)}}\sum_{n=1}^\infty e^{-b^2n^2/(\sigma + it)} \\
        &+ \frac{ab}{ \pi(\sigma + it)} \sum_{m,n \geq 1} e^{-(a^2m^2 +b^2n^2)/(\sigma + it)}. 
    \end{split}
\end{equation}
We use the following convention of the Fourier transform:
\begin{equation} \label{fourier transform convention}
    \widehat f(x) = \int_{-\infty}^\infty f(t) e^{-i t x} dt.
\end{equation}
For $E,\gamma > 0$, denote by $\varphi_{E,\gamma}(t) = \varphi(t)$ the entire function which is the inverse Fourier transform of the function 
\begin{equation*}
    x \mapsto (|x|- E)_-^\gamma = 
    \begin{cases}
        (E - |x|)^\gamma, \quad &|x| \leq E, \\
        0, \quad &|x| > E.
    \end{cases}
\end{equation*}
%\red{Move lemma to appendix?}
\begin{lemma} \label{lemma varphi estimate}
    We have
    \begin{equation} \label{eq: varphi estimate}
        \varphi(t) = O_{E,\gamma}((1 + |t|)^{-\min(\gamma + 1,2)}), \quad t \in \R.
    \end{equation}
    In particular, $\varphi \in L^1(\R)$ for every $E,\gamma > 0$.
\end{lemma}
\begin{proof}
    We have
    \begin{equation*}
        \varphi(t) = \frac{1}{2\pi} \int_{-\infty}^\infty  (|x|- E)_-^\gamma e^{ixt} dx = \frac{1}{\pi} \int_0^E (E - x)^\gamma \cos(xt) dx.
    \end{equation*}
    For bounded $t$, there is nothing to prove, since
    \begin{equation*}
        |\varphi(t)| \leq \frac{1}{\pi} \int_0^E (E - x)^\gamma dx = \frac{E^{\gamma + 1}}{\pi(\gamma + 1)}.
    \end{equation*}
    It is therefore enough to show that
    \begin{equation*}
        \varphi(t) = O_{E,\gamma}(|t|^{-\min(\gamma + 1,2)}), \quad |t| \to \infty.
    \end{equation*}
    Integration by parts gives
    \begin{equation*}
        \varphi(t) = \frac{\gamma}{\pi t} \int_0^E (E - x)^{\gamma - 1} \sin(tx) dx,
    \end{equation*}
    so it is enough to show that
    \begin{equation*}
        \int_0^E (E - x)^{\gamma - 1} \sin(tx) dx = O_{E,\gamma}(|t|^{-\min(\gamma, 1)}), \quad |t| \to \infty.
    \end{equation*}
    If $\gamma = 1$, we directly compute
    \begin{equation*}
        \int_0^E \sin(tx) dx = \frac{1 - \cos(Et)}{t} = O(|t|^{-1}), \quad |t| \to \infty.
    \end{equation*}
    Now suppose $\gamma > 1$. Integration by parts gives
    \begin{equation*}
        \int_0^E (E - x)^{\gamma - 1} \sin(tx) dx = \frac{E^{\gamma - 1}}{t} - \frac{\gamma - 1}{t} \int_0^E (E - x)^{\gamma - 2} \cos(tx) dx.
    \end{equation*}
    The integral is uniformly bounded since
    \begin{equation*}
        \bigg|\int_0^E (E - x)^{\gamma - 2} \cos(tx) dx\bigg| \leq \int_0^E (E - x)^{\gamma - 2} dx = \frac{E^{\gamma - 1}}{\gamma - 1},
    \end{equation*}
    hence the lemma follows for $\gamma > 1$. Finally, suppose $0 < \gamma < 1$. We make the substitutions $u = E - x$ and $v = |t|u$ to obtain
    \begin{equation*}
        \int_0^E (E - x)^{\gamma - 1} \sin(tx) dx = \int_0^E u^{\gamma - 1} \sin(Et - tu) du = |t|^{-\gamma} \int_0^{E|t|} v^{\gamma - 1} \sin(Et - \text{sgn}(t)v) dv.
    \end{equation*}
    Thus, it is enough to show that
    %\begin{equation*}
        %\int_0^{R} v^{\gamma - 1} \sin(Et - \text{sgn}(t)v) dv
    %\end{equation*}
    %is uniformly bounded in $R > 0$. Moreover, since 
    %\begin{align*}
        %\int_0^{R} v^{\gamma - 1} \sin(Et - \text{sgn}(t)v) dv = \text{Im}\bigg(e^{i E t}\int_0^{R} v^{\gamma - 1} e^{-i \text{sgn}(t)v} dv\bigg),
    %\end{align*}
    %it is enough to show that 
    \begin{equation*}
        \int_0^R v^{\gamma - 1} e^{\pm iv} dv
    \end{equation*}
    is uniformly bounded in $R > 0$. In fact, since
    \begin{equation*}
        \bigg|\int_0^1 v^{\gamma - 1} e^{\pm iv} dv\bigg| \leq \int_0^1 v^{\gamma - 1} dv = \frac{1}{\gamma},
    \end{equation*}
    it is enough to show that
    \begin{equation*}
        \int_1^R v^{\gamma - 1} e^{\pm iv} dv
    \end{equation*}
    is uniformly bounded in $R > 1$. Integration by parts gives
    \begin{equation*}
        \int_1^R v^{\gamma - 1} e^{\pm iv} dv = \frac{R^{\gamma-1}e^{\pm i R} - e^{\pm i}}{\pm i} \mp \frac{\gamma - 1}{i} \int_1^R v^{\gamma - 2}e^{\pm iv} dv.
    \end{equation*}
    Since $0 < \gamma < 1$, it follows that
    \begin{equation*}
        \bigg|\frac{R^{\gamma-1}e^{\pm i R} - e^{\pm i}}{\pm i}\bigg| \leq 2,
    \end{equation*}
    and 
    \begin{equation*}
        \bigg|\int_1^R v^{\gamma - 2}e^{\pm iv} dv\bigg| \leq \int_1^R v^{\gamma - 2} dv = \frac{1 - R^{\gamma - 1}}{1 - \gamma} \leq \frac{1}{1 - \gamma}.
    \end{equation*}
    This completes the proof of the lemma.
\end{proof}

Now multiply $\sum_{j = 0}^\infty e^{-(\sigma + it) \lambda_j}$ by $\varphi(t)$ and integrate over $t \in \R$. We get 
\begin{equation*}
    \begin{split}
        \int_{-\infty}^\infty \varphi(t) \sum_{j = 0}^\infty e^{-(\sigma + it) \lambda_j} dt &= \sum_{j = 0}^\infty \int_{-\infty}^\infty  \varphi(t) e^{-(\sigma + it) \lambda_j} dt  \\
        &=  \sum_{j=0}^\infty e^{-\sigma \lambda_j}  \mathcal{F} \varphi(\lambda_j) \\
        &= \sum_{\lambda_j \leq E} e^{-\sigma \lambda_j} (E - \lambda_j)^\gamma.
    \end{split}
\end{equation*}
Since $\varphi \in L^1(\R)$ and $\sum_j e^{-\lambda_j \sigma} < \infty$, the interchange of the sum and the integral is justified by the dominated convergence theorem. Thus, by \eqref{complex heat trace id},
\begin{equation} \label{sigma frozen}
    \begin{split}
         \sum_{\lambda_j \leq E} e^{-\sigma \lambda_j} (E - \lambda_j)^\gamma &= \frac{ab}{4\pi} \int_{-\infty}^\infty  \frac{\varphi(t)}{\sigma + it} dt \pm \frac{a + b}{4 \sqrt{\pi}} \int_{-\infty}^\infty  \frac{\varphi(t)}{\sqrt{\sigma + it}} dt + \frac{1}{4} \int_{-\infty}^\infty  \varphi(t) dt \\
         &+ \int_{-\infty}^\infty  \varphi(t) \mathcal{R}_\square(\sigma + it) dt. 
     \end{split}
\end{equation}
We now analyze the limit as $\sigma \to 0^+$ term by term in \eqref{sigma frozen}. For the left-hand side, it is clear that 
\begin{equation*}
    \lim_{\sigma \to 0^+} \sum_{\lambda_j \leq E} e^{-\sigma \lambda_j} (E - \lambda_j)^\gamma = \sum_{\lambda_j \leq E} (E - \lambda_j)^\gamma = R_\square^\gamma(E).
\end{equation*}
For the right-hand side of \eqref{sigma frozen},
\begin{equation*}
    \int_{-\infty}^\infty \varphi(t) dt = \mathcal{F} \varphi(0) = E^\gamma,
\end{equation*}
and for $m \geq 1$ we compute 
\begin{align*}
    \lim_{\sigma \to 0^+}\int_{-\infty}^\infty \frac{\varphi(t)}{(\sigma + it)^{m/2}} dt &= \lim_{\sigma \to 0^+}\int_{-\infty}^\infty \varphi(t) \frac{1}{\Gamma(m/2)} \int_0^\infty x^{m/2-1} e^{-x(\sigma + it)}dxdt \\
    &= \lim_{\sigma \to 0^+} \frac{1}{\Gamma(m/2)} \int_{0}^\infty x^{m/2-1}e^{-\sigma x} \int_{-\infty}^\infty \varphi(t)e^{-i x t} dtdx \\
    &= \lim_{\sigma \to 0^+} \frac{1}{\Gamma(m/2)} \int_{0}^E x^{m/2-1} e^{-\sigma x}(E - x)^\gamma dx \\
    &= \frac{1}{\Gamma(m/2)} \int_{0}^E x^{m/2-1} (E - x)^\gamma dx \\
    &= \frac{1}{\Gamma(m/2)}E^{m/2+\gamma} \int_0^1 x^{m/2-1}(1-x)^\gamma dx \\
    &= \frac{\Gamma(\gamma + 1)}{\Gamma(m/2 + \gamma + 1)}E^{m/2+\gamma}.
\end{align*}
Note that for every $\sigma > 0$, Fubini's theorem applies since
\begin{equation*}
    \int_{-\infty}^\infty \int_0^\infty |\varphi(t)x^{m/2-1}e^{-x(\sigma + it)}| dxdt = ||\varphi||_{L^1(\R)} \int_0^\infty x^{m/2-1}e^{-\sigma x} dx < \infty.
\end{equation*}
Moreover, %since $\mathcal{F}\varphi(x)$ has compact support, 
the interchange of the limit and the integral is justified by the dominated convergence theorem, since $x^{m/2-1}(E-x)^\gamma \in L^1(0,E)$ for $m \geq 1$. Thus, by sending $\sigma \to 0^+$ in \eqref{sigma frozen}, we find that the Riesz means of the rectangle of side lengths $a$ and $b$ are
\begin{equation} \label{find asymptotic for}
    R_\square^\gamma(E) = \frac{ab}{4\pi(\gamma + 1)} E^{\gamma + 1} \pm \frac{a + b}{4\sqrt{\pi}}\frac{\Gamma(\gamma + 1)}{\Gamma(\gamma + 3/2)} E^{\gamma + 1/2} + \frac{1}{4} E^\gamma + \lim_{\sigma \to 0^+} \int_{-\infty}^\infty \varphi(t) \mathcal{R}_\square(\sigma + it) dt.
\end{equation}
It remains to analyze the remainder
\begin{equation*}
    \lim_{\sigma \to 0^+} \int_{-\infty}^\infty \varphi(t) \mathcal{R}_\square(\sigma + it) dt.
\end{equation*}
We start by bounding the partial sums appearing in $\mathcal{R}_\square(\sigma + it)$, which will allow us to justify certain interchanges between integrals and infinite sums. 
\begin{lemma} \label{lemma: sharp exponential bound on partial sums}
    For $\sigma > 0$, we have
    \begin{equation} \label{eq: sharp exponential bound partial sums}
        \sum_{n = 1}^N e^{-n^2/(\sigma + it)}  = O_\sigma((1 + |t|)^{1/2}),
    \end{equation}
    uniformly in $N \in \NN$ and $t \in \R$.
\end{lemma}
\begin{proof}
    Write
    \begin{equation*}
        \alpha \defeq \frac{\sigma}{\sigma^2 + t^2}, \quad \beta \defeq \frac{t}{\sigma^2 + t^2},
    \end{equation*}
    so that
    \begin{equation*}
        e^{-n^2/(\sigma + it)} = e^{-\alpha n^2} e^{i \beta n^2}.
    \end{equation*}
    For bounded $t$, say $|t| \leq T$, we have $\alpha \geq \sigma/(\sigma^2 + T^2)$, hence
    \begin{equation*}
        \bigg|\sum_{n=1}^N e^{-n^2/(\sigma+it)}\bigg| \leq \sum_{n=1}^\infty e^{-\alpha n^2} \leq 
        \sum_{n=1}^\infty e^{- n^2 \sigma/(\sigma^2+T^2)} \lesssim_{\sigma} 1.
    \end{equation*}
    Thus, it is enough to show that 
    \begin{equation*}
        \sum_{n = 1}^N e^{-n^2/(\sigma + it)}  = O_\sigma(|t|^{1/2}), \quad |t| \to \infty,
    \end{equation*}
    uniformly in $N \in \NN$. To this end, we use summation by parts to compute
    \begin{equation*}
        \sum_{n = 1}^N e^{-n^2/(\sigma + it)} = \sum_{n = 1}^N e^{-\alpha n^2} e^{i \beta n^2} = e^{-\alpha N^2}A_N(\beta) + \sum_{n = 1}^{N-1} (e^{-\alpha n^2} - e^{-\alpha(n+1)^2})A_n(\beta),
    \end{equation*}
    where $A_n(\beta) \defeq \sum_{k = 1}^n e^{i \beta k^2}$. By the van der Corput estimate \cite[Theorem 2.2]{grko91},
    \begin{equation*}
        A_n(\beta) \lesssim n |\beta|^{1/2} + |\beta|^{-1/2}.
    \end{equation*}
    Hence
    \begin{align*}
        \bigg|\sum_{n = 1}^N e^{-n^2/(\sigma + it)}\bigg| &\lesssim e^{-\alpha N^2}(N |\beta|^{1/2} + |\beta|^{-1/2}) + \sum_{n = 1}^{N-1} (e^{-\alpha n^2} - e^{-\alpha(n+1)^2})(n |\beta|^{1/2} + |\beta|^{-1/2}) \\
        &= |\beta|^{1/2} \sum_{n = 1}^N e^{-\alpha n^2} + |\beta|^{-1/2} e^{-\alpha} \\
        &\lesssim \sqrt{\frac{|\beta|}{\alpha}} + \sqrt{\frac{1}{|\beta|}} \\
        &= \sqrt{\frac{|t|}{\sigma}} + \sqrt{\frac{\sigma^2}{|t|} + |t|} \\
        &\lesssim_{\sigma} |t|^{1/2}, \quad |t| \to \infty.
    \end{align*}
    This proves the lemma.
\end{proof}
Now, consider the single sums in \eqref{R of sigma plus it} involving $(\sigma + it)^{-1/2}$. By Lemmas \ref{lemma varphi estimate} and \ref{lemma: sharp exponential bound on partial sums},
\begin{equation} \label{eq: varphi sum estimate} 
    \bigg|\varphi(t) \frac{1}{\sqrt{\sigma + it}} \sum_{m = 1}^M e^{-a^2m^2/(\sigma + it)}\bigg| \lesssim_\sigma(1 + |t|)^{-\min(\gamma+1,2)}, \quad M \in \NN, \quad t \in \R.
\end{equation}
In particular, for every fixed $\sigma,\gamma > 0$, the right-hand side of \eqref{eq: varphi sum estimate} is in $L^1(\R)$. Thus, the dominated convergence theorem yields
\begin{equation} \label{eq: DCT -1/2}
    \begin{split}
        \int_{-\infty}^\infty \varphi(t) \frac{a}{2\sqrt{\pi(\sigma + it)}} \sum_{m=1}^\infty e^{-a^2m^2/(\sigma + it)} dt = \frac{a}{2\sqrt{\pi}} \sum_{m = 1}^\infty \int_{-\infty}^\infty \varphi(t) \frac{1}{\sqrt{\sigma + it}} e^{-a^2m^2/(\sigma + it)} dt.
    \end{split}
\end{equation}
By \cite[3.896.4]{grry07} and the substitution $x = s^2$, we have
\begin{equation*}
    \frac{1}{\sqrt{z}}e^{-y/z} = \frac{1}{\sqrt{\pi}} \int_0^\infty x^{-1/2} e^{-zx} \cos(2\sqrt{xy}) dx, \quad \text{Re}(z) > 0, \quad y > 0.
\end{equation*}
This gives 
\begin{equation} \label{eq single sum -1/2 step 1}
    \begin{split}
        &\frac{a}{2\sqrt{\pi}} \sum_{m = 1}^\infty \int_{-\infty}^\infty \varphi(t) \frac{1}{\sqrt{\sigma + it}} e^{-a^2m^2/(\sigma + it)} dt \\
        &= \frac{a}{2\pi} \sum_{m = 1}^\infty\int_{-\infty}^\infty \varphi(t) \int_0^\infty x^{-1/2} e^{-(\sigma + it)x} \cos(2am\sqrt{x}) dx dt \\
        &= \frac{a}{2\pi} \sum_{m = 1}^\infty \int_0^\infty x^{-1/2} e^{-\sigma x} \cos(2am\sqrt{x}) \int_{-\infty}^\infty \varphi(t) e^{-i t x} dt dx \\
        &= \frac{a}{2\pi} \sum_{m = 1}^\infty \int_0^E (E-x)^\gamma x^{-1/2} e^{-\sigma x} \cos(2am\sqrt{x}) dx \\
        &= \frac{a}{\pi} \sum_{m = 1}^\infty \int_0^{\sqrt{E}} (E-x^2)^\gamma e^{-\sigma x^2} \cos(2amx) dx.
    \end{split}
\end{equation}
In the second equality, we interchange the integrals. This is justified by Fubini's theorem, since 
\begin{equation*}
    \int_{-\infty}^\infty \int_0^\infty |\varphi(t) x^{-1/2} e^{-(\sigma + it)x} \cos(2 a m \sqrt{x})| dxdt \leq ||\varphi||_{L^1(\R)} \int_0^\infty x^{-1/2} e^{-\sigma x} dx 
    < \infty.
\end{equation*}
Now we want to send $\sigma \to 0^+$ in \eqref{eq single sum -1/2 step 1}. In Appendix \ref{estimate appendix new}, we show that (see \eqref{cosine decay})
\begin{equation} \label{eq: integral cosine estimate}
    \bigg|\int_0^{\sqrt E} (E-x^2)^\gamma e^{-\sigma x^2} \cos(2amx)\,dx\bigg| \leq Cm^{-\gamma - 1},
\end{equation}
for some $C > 0$ independent of $m \in \NN$ and $0 \leq \sigma \leq 1$. Since $m^{-\gamma - 1}$ is summable, the dominated convergence theorem allows us to interchange the limit with the infinite sum:
\begin{equation*} 
    \begin{split}
        \lim_{\sigma \to 0^+} \frac{a}{\pi} \sum_{m = 1}^\infty \int_0^{\sqrt{E}} (E-x^2)^\gamma e^{-\sigma x^2} \cos(2amx) dx \\
        = \frac{a}{\pi} \sum_{m = 1}^\infty \lim_{\sigma \to 0^+} \int_0^{\sqrt{E}} (E-x^2)^\gamma e^{-\sigma x^2} \cos(2amx) dx.
    \end{split}
\end{equation*}
Then, another application of the dominated convergence theorem justifies interchanging the limit with the integral, since the integrand is dominated by $(E-x^2)^\gamma  \in L^1(0,\sqrt{E})$. Thus, we obtain
\begin{equation} \label{eq single sum -1/2}
    \begin{split}
        &\lim_{\sigma \to 0^+} \int_{-\infty}^\infty \varphi(t) \frac{a}{2\sqrt{\pi(\sigma + it)}} \sum_{m=1}^\infty e^{-a^2m^2/(\sigma + it)} dt \\
        &= \frac{a}{\pi} \sum_{m = 1}^\infty \int_0^{\sqrt{E}} (E-x^2)^\gamma  \cos(2amx) dx \\
        &= \frac{\Gamma(\gamma + 1)}{2a^{\gamma-1/2}\sqrt{\pi}} E^{\gamma/2 + 1/4} \sum_{m = 1}^\infty \frac{J_{\gamma + 1/2}(2am\sqrt{E})}{m^{\gamma+1/2}}.
    \end{split}
\end{equation}
%\begin{equation} \label{eq single sum -1/2}
    %\begin{split}
        %&\frac{a}{2\sqrt{\pi}} \lim_{\sigma \to 0^+} \sum_{m = 1}^\infty \int_{-\infty}^\infty \varphi(t) \frac{1}{\sqrt{\sigma + it}} e^{-a^2m^2/(\sigma + it)} dt \\
        %&= \frac{a}{2\pi} \lim_{\sigma \to 0^+} \int_{-\infty}^\infty \varphi(t) \sum_{m = 1}^\infty \int_0^\infty x^{-1/2} e^{-(\sigma + it)x} \cos(2am\sqrt{x}) dx dt \\
        %&= \frac{a}{2\pi} \lim_{\sigma \to 0^+} \sum_{m = 1}^\infty\int_{-\infty}^\infty \varphi(t) \int_0^\infty x^{-1/2} e^{-(\sigma + it)x} \cos(2am\sqrt{x}) dx dt \\
        %&= \frac{a}{2\pi} \lim_{\sigma \to 0^+} \sum_{m = 1}^\infty \int_0^\infty x^{-1/2} e^{-\sigma x} \cos(2am\sqrt{x}) \int_{-\infty}^\infty \varphi(t) e^{-i t x} dt dx \\
        %&= \frac{a}{2\pi} \lim_{\sigma \to 0^+} \sum_{m = 1}^\infty \int_0^E (E-x)^\gamma x^{-1/2}e^{-\sigma x} \cos(2am\sqrt{x}) dx \\
        %&= \frac{a}{2\pi} \sum_{m = 1}^\infty \lim_{\sigma \to 0^+} \int_0^E (E-x)^\gamma x^{-1/2}e^{-\sigma x} \cos(2am\sqrt{x}) dx \\
        %&= \frac{a}{2\pi} \sum_{m = 1}^\infty \int_0^E (E-x)^\gamma x^{-1/2} \cos(2am\sqrt{x}) dx \\
        %&= \frac{a}{\pi} \sum_{m = 1}^\infty \int_0^{\sqrt{E}} (E-x^2)^\gamma \cos(2amx) dx \\
        %&= \frac{\Gamma(\gamma + 1)}{2a^{\gamma-1/2}\sqrt{\pi}} E^{\gamma/2 + 1/4} \sum_{m = 1}^\infty \frac{J_{\gamma + 1/2}(2am\sqrt{E})}{m^{\gamma+1/2}}.
    %\end{split}
%\end{equation}
Here, $J_\nu(z)$ denotes the Bessel function of the first kind of order $\nu \ge 0$ \cite[\S 10.2(ii)]{dlmf}. The last equality in \eqref{eq single sum -1/2} follows from %the standard integral identity for $J_\nu$ 
\cite[\S 10.9(i), equation 10.9.4]{dlmf}. A similar calculation yields
\begin{equation} \label{eq2 single sum -1/2}
    \lim_{\sigma \to 0^+} \int_{-\infty}^\infty \varphi(t) \frac{b}{2\sqrt{\pi(\sigma + it)}} \sum_{n=1}^\infty e^{-b^2n^2/(\sigma + it)} dt = \frac{\Gamma(\gamma + 1)}{2b^{\gamma-1/2}\sqrt{\pi}} E^{\gamma/2 + 1/4} \sum_{n = 1}^\infty \frac{J_{\gamma + 1/2}(2bn\sqrt{E})}{n^{\gamma+1/2}}.
\end{equation}
By the asymptotic expansion \cite[\S 10.17(i)]{dlmf}
\begin{equation} \label{eq: J_2 asymptotic}
    J_\nu(x) = \sqrt{\frac{2}{\pi x}}\cos\bigg(x - \frac{(2\nu + 1)\pi}{4}\bigg) + O(x^{-3/2}), \quad x \to \infty,
\end{equation}
the summands on the right-hand sides of \eqref{eq single sum -1/2} and \eqref{eq2 single sum -1/2} are $O(m^{-\gamma - 1})$ and $O(n^{-\gamma - 1})$, respectively. Thus, both series converge absolutely for every $\gamma > 0$.

Next, consider the single sums in \eqref{R of sigma plus it} involving $(\sigma + it)^{-1}$. By Lemmas \ref{lemma varphi estimate} and \ref{lemma: sharp exponential bound on partial sums},
\begin{equation} \label{eq: varphi sum estimate (sigma it) -1}
    \bigg|\varphi(t) \frac{1}{\sigma + it} \sum_{m = 1}^M e^{-a^2m^2/(\sigma + it)}\bigg| = O_\sigma((1 + |t|)^{-\min(\gamma + \frac{3}{2},\frac{5}{2})}), \quad M \in \NN, \quad t \in \R,
\end{equation}
and the dominated convergence theorem therefore gives
\begin{equation} \label{eq: DCT -1}
    \begin{split}
    \int_{-\infty}^\infty \varphi(t) \frac{ab}{2\pi(\sigma + it)} \sum_{m=1}^\infty e^{-a^2m^2/(\sigma + it)} dt = \frac{ab}{2\pi} \sum_{m=1}^\infty\int_{-\infty}^\infty \varphi(t) \frac{1}{\sigma + it} e^{-a^2m^2/(\sigma + it)} dt
    \end{split}
\end{equation}
for every fixed $\sigma > 0$. We then use
\begin{equation} \label{eq: Bessel identity}
    \frac{1}{z}e^{-y/z} = \int_0^\infty e^{-zx} J_0(2\sqrt{xy}) dx, \quad \text{Re}(z) > 0,
\end{equation}
obtained as a special case of \cite[6.615]{grry07}, to compute 
\begin{equation} \label{eq single sum -1 step1}
    \begin{split}
        &\frac{ab}{2\pi} \sum_{m=1}^\infty\int_{-\infty}^\infty \varphi(t) \frac{1}{\sigma + it} e^{-a^2m^2/(\sigma + it)} dt \\
        &= \frac{ab}{2\pi} \sum_{m = 1}^\infty \int_{-\infty}^\infty \int_0^\infty \varphi(t) e^{-(\sigma + it)x} J_0(2am\sqrt{x}) dxdt \\
        &= \frac{ab}{2\pi} \sum_{m = 1}^\infty \int_0^\infty e^{-\sigma x} J_0(2 a m \sqrt{x}) \int_{-\infty}^\infty \varphi(t) e^{-itx} dt dx \\
        &= \frac{ab}{2\pi} \sum_{m = 1}^\infty \int_0^E (E - x)^\gamma e^{-\sigma x} J_0(2 a m \sqrt{x}) dx \\
        &= \frac{ab}{\pi} \sum_{m = 1}^\infty \int_0^{\sqrt{E}} x(E - x^2)^\gamma e^{-\sigma x^2} J_0(2 a m x) dx.
    \end{split}
\end{equation}
The interchange of the integrals is again justified by Fubini's theorem. By \eqref{bessel bound}, we can as before use the dominated convergence theorem twice to interchange the limit $\sigma \to 0^+$ with the sum and integral. This yields
\begin{equation} \label{eq single sum -1}
    \begin{split}
        \lim_{\sigma \to 0^+}\int_{-\infty}^\infty \varphi(t) \frac{ab}{2\pi(\sigma + it)} \sum_{m=1}^\infty e^{-a^2m^2/(\sigma + it)} dt &= \frac{ab}{\pi} \sum_{m = 1}^\infty \int_0^{\sqrt{E}} x(E - x^2)^\gamma J_0(2 a m x) dx \\
        &= \frac{b\Gamma(\gamma + 1)}{2\pi a^\gamma} E^{(\gamma + 1)/2} \sum_{m = 1}^\infty \frac{J_{\gamma + 1}(2am\sqrt{E})}{m^{\gamma + 1}}.
    \end{split}
\end{equation}
%\begin{equation} \label{eq single sum -1}
    %\begin{split}
        %&\lim_{\sigma \to 0^+} \int_{-\infty}^\infty \varphi(t) \frac{ab}{2\pi(\sigma + it)} \sum_{m=1}^\infty e^{-a^2m^2/(\sigma + it)} dt \\
        %&\frac{ab}{2\pi}\lim_{\sigma \to 0^+} \sum_{m=1}^\infty\int_{-\infty}^\infty \varphi(t) \frac{1}{\sigma + it} e^{-a^2m^2/(\sigma + it)} dt \\
        %&= \frac{ab}{2\pi}\lim_{\sigma \to 0^+} \sum_{m = 1}^\infty \int_{-\infty}^\infty \int_0^\infty \varphi(t) e^{-(\sigma + it)x} J_0(2am\sqrt{x}) dxdt \\
        %&= \frac{ab}{2\pi}\lim_{\sigma \to 0^+} \sum_{m = 1}^\infty \int_0^\infty e^{-\sigma x} J_0(2 a m \sqrt{x}) \int_{-\infty}^\infty \varphi(t) e^{-itx} dt dx \\
        %&= \frac{ab}{2\pi}\lim_{\sigma \to 0^+} \sum_{m = 1}^\infty \int_0^E e^{-\sigma x} (E - x)^\gamma J_0(2 a m \sqrt{x}) dx \\
        %&= \frac{ab}{2\pi} \sum_{m = 1}^\infty \lim_{\sigma \to 0^+} \int_0^E e^{-\sigma x} (E - x)^\gamma J_0(2 a m \sqrt{x}) dx \\
        %&= \frac{ab}{2\pi} \sum_{m = 1}^\infty \int_0^E (E - x)^\gamma J_0(2 a m \sqrt{x}) dx \\
        %&= \frac{ab}{\pi} \sum_{m = 1}^\infty \int_0^{\sqrt{E}} x(E - x^2)^\gamma J_0(2 a m x) dx \\
        %&= \frac{b\Gamma(\gamma + 1)}{2\pi a^\gamma} E^{(\gamma + 1)/2} \sum_{m = 1}^\infty \frac{J_{\gamma + 1}(2am\sqrt{E})}{m^{\gamma + 1}}.
    %\end{split}
%\end{equation}
%The interchange of the limit and the infinite sum is again justified in Appendix \ref{estimate appendix}. As before, the interchanges of the integrals as well as the limit with the integral follow from Fubini's theorem and the dominated convergence theorem, respectively. 
The last equality uses Sonine's integral identity \cite{we}. Repeating this argument,
\begin{equation} \label{eq2 single sum -1}
    \lim_{\sigma \to 0^+} \int_{-\infty}^\infty \varphi(t) \frac{ab}{2\pi(\sigma + it)} \sum_{n=1}^\infty e^{-b^2n^2/(\sigma + it)} dt = \frac{a\Gamma(\gamma + 1)}{2\pi b^\gamma}E^{(\gamma + 1)/2} \sum_{n = 1}^\infty \frac{J_{\gamma+1}(2 b n \sqrt{E})}{n^{\gamma + 1}}.
\end{equation}
By \eqref{eq: J_2 asymptotic}, the series on the right-hand sides of \eqref{eq single sum -1} and \eqref{eq2 single sum -1} converge absolutely for every $\gamma > 0$. 

Finally, consider the double sum in \eqref{R of sigma plus it}. By Lemmas \ref{lemma varphi estimate} and \ref{lemma: sharp exponential bound on partial sums},
\begin{equation*}
    \bigg| \varphi(t) \frac{1}{\sigma + it} \sum_{1 \leq m,n \leq N} e^{-(a^2m^2 + b^2n^2)/(\sigma + it)} \bigg| = O_\sigma((1 + |t|)^{-\min(\gamma + 1, 2)}),
\end{equation*}
uniformly in $N \in \NN$ and $t \in \R$. Since the right-hand side is integrable, the dominated convergence theorem allows us to interchange the limit with the integral. Thus, with the notation
\begin{equation*}
    \sum_{m,n \geq 1}^\square := \lim_{N \to \infty} \sum_{1 \leq m,n \leq N},
\end{equation*}
%keeping the rectangular truncation explicit, 
we obtain
\begin{equation} \label{eq: DCT double}
    \begin{split}
        &\int_{-\infty}^\infty \varphi(t) \frac{ab}{\pi(\sigma + it)} \sum_{m,n \geq 1} e^{-(a^2m^2 + b^2n^2)/(\sigma + it)} dt \\
        %&=\frac{ab}{\pi}\int_{-\infty}^\infty \varphi(t) \frac{1}{\sigma + it} \lim_{N \to \infty} \sum_{1 \leq m,n \leq N} e^{-(a^2m^2 + b^2n^2)/(\sigma + it)} dt \\
        &= \frac{ab}{\pi} \sum_{m,n \geq 1}^\square \int_{-\infty}^\infty \frac{\varphi(t)}{\sigma + it} e^{-(a^2m^2+b^2n^2)/(\sigma + it)} dt \\
        &= \frac{ab}{\pi} \sum_{m,n \geq 1}^\square \int_{-\infty}^\infty \varphi(t) \int_0^\infty e^{-(\sigma + it)x} J_0\big(2\sqrt{x(a^2m^2 + b^2n^2)}\big) dxdt \\
        &= \frac{ab}{\pi} \sum_{m,n \geq 1}^\square \int_{0}^\infty e^{-\sigma x} J_0\big(2\sqrt{x(a^2m^2 + b^2n^2)}\big) \int_{-\infty}^\infty \varphi(t) e^{-itx}  dtdx \\
        &= \frac{ab}{\pi} \sum_{m,n \geq 1}^\square \int_{0}^E (E-x)^\gamma e^{-\sigma x} J_0\big(2\sqrt{x(a^2m^2 + b^2n^2)}\big) dx \\
        &= \frac{2ab}{\pi} \sum_{m,n \geq 1}^\square \int_0^{\sqrt{E}} x (E - x^2)^\gamma e^{-\sigma x^2} J_0\big(2x\sqrt{a^2m^2 + b^2n^2}\big) dx,
    \end{split}
\end{equation}
for fixed $\sigma > 0$. However, unlike in the previous cases, the series on the right-hand side need not converge absolutely. Indeed, by \eqref{eq: J_sigma first term}, we have
\begin{equation} \label{eq: integral J_0 estimate}
    \begin{split}
        &\int_0^{\sqrt{E}} x (E - x^2)^\gamma e^{-\sigma x^2} J_0\big(2x\sqrt{a^2m^2 + b^2n^2}\big) dx \\
        &= e^{-\sigma E}C_{E,\gamma} \frac{\cos\big(2\sqrt{E(a^2m^2+b^2n^2)} - \frac{(2\gamma + 3)\pi}{4}\big)}{(a^2m^2 + b^2n^2)^{\gamma/2 + 3/4}} + O((a^2m^2+b^2n^2)^{-\gamma/2-1}),
    \end{split}
\end{equation}
uniformly in $0 \leq \sigma \leq 1$. Since the remainder is absolutely summable over $m,n \geq 1$, it follows that the right-hand side of \eqref{eq: DCT double} is absolutely convergent if and only if 
\begin{equation*}
    \frac{\gamma}{2} + \frac{3}{4} > 1 \iff \gamma > 1/2.
\end{equation*}
For $0 < \gamma \leq 1/2$, the rectangular convergence in \eqref{eq: DCT double} is therefore only conditional. We next use this convergence to show that the corresponding Bessel series also converges under rectangular summation. Write
\begin{equation} \label{eq: h_sigma calculation}
    \begin{split}
        &\int_0^{\sqrt{E}} x (E - x^2)^\gamma e^{-\sigma x^2} J_0\big(2x\sqrt{a^2m^2 + b^2n^2}\big) dx \\
        &= e^{-\sigma E}\int_0^{\sqrt{E}} x (E - x^2)^\gamma e^{\sigma (E - x^2)} J_0\big(2x\sqrt{a^2m^2 + b^2n^2}\big) dx \\
        &= e^{-\sigma E}\int_0^{\sqrt{E}} x (E - x^2)^\gamma J_0\big(2x\sqrt{a^2m^2 + b^2n^2}\big) dx \\
        &+ e^{-\sigma E}\int_0^{\sqrt{E}} x (E - x^2)^\gamma (e^{\sigma (E - x^2)} - 1) J_0\big(2x\sqrt{a^2m^2 + b^2n^2}\big) dx \\
        &= \frac{1}{2}e^{-\sigma E} \, \Gamma(\gamma + 1) E^{(\gamma + 1)/2} \frac{J_{\gamma + 1}\big(2\sqrt{E(a^2m^2 + b^2n^2)}\big)}{(a^2m^2 + b^2n^2)^{(\gamma + 1)/2}} \\
        &+ \sigma e^{-\sigma E} \int_0^{\sqrt{E}} x (E - x^2)^{\gamma+1} h_\sigma(E-x^2) J_0\big(2x\sqrt{a^2m^2 + b^2n^2}\big) dx.
    \end{split}
\end{equation}
Here, $h_\sigma(y) = \int_0^1 e^{s\sigma y} ds$. By \eqref{bessel bound}, the integral appearing in the second term on the right-hand side of \eqref{eq: h_sigma calculation} satisfies
\begin{equation} \label{eq: integral h_sigma estimate}
    \int_0^{\sqrt{E}} x (E - x^2)^{\gamma+1} h_\sigma(E-x^2) J_0\big(2x\sqrt{a^2m^2 + b^2n^2}\big) dx = O\big((a^2m^2 + b^2n^2)^{-\gamma/2-5/4}\big),
\end{equation}
uniformly in $0 \leq \sigma \leq 1$, and is therefore absolutely summable over $m,n \geq 1$. On the other hand, for every $\sigma > 0$, the left-hand side of \eqref{eq: h_sigma calculation} is summable under rectangular summation by \eqref{eq: DCT double}. It follows that
\begin{equation*}
    \sum_{m,n \geq 1}^\square \frac{J_{\gamma + 1}\big(2\sqrt{E(a^2m^2 + b^2n^2)}\big)}{(a^2m^2 + b^2n^2)^{(\gamma + 1)/2}}
\end{equation*}
exists for every $\gamma > 0$. Finally, letting $\sigma \to 0^+$ in \eqref{eq: DCT double} and using \eqref{eq: h_sigma calculation} gives
\begin{equation} \label{eq: limit double sum final}
    \begin{split}
        &\lim_{\sigma \to 0^+} \int_{-\infty}^\infty \varphi(t) \frac{ab}{\pi(\sigma + it)} \sum_{m,n \geq 1} e^{-(a^2m^2 + b^2n^2)/(\sigma + it)} dt \\
        &= \frac{ab \Gamma(\gamma + 1)}{\pi} E^{(\gamma + 1)/2} \sum_{m,n \geq 1}^\square \frac{J_{\gamma + 1}\big(2\sqrt{E(a^2m^2 + b^2n^2)}\big)}{(a^2m^2 + b^2n^2)^{(\gamma + 1)/2}},
    \end{split}
\end{equation}
for every $\gamma > 0$. %If $\gamma > 1/2$, the right-hand side of \eqref{eq: limit double sum final} converges absolutely, and we simply obtain
%\begin{equation} \label{eq double sum}
    %\begin{split}
        %&\lim_{\sigma \to 0^+} \int_{-\infty}^\infty \varphi(t) \frac{ab}{\pi(\sigma + it)} \sum_{m,n \geq 1} e^{-(a^2m^2 + b^2n^2)/(\sigma + it)} dt \\
        %&= \frac{ab\Gamma(\gamma + 1)}{\pi}E^{(\gamma+1)/2} \sum_{m,n \geq 1} \frac{J_{\gamma + 1}\big(2\sqrt{E(a^2m^2 + b^2n^2)}\big)}{(a^2m^2 + b^2n^2)^{(\gamma + 1)/2}}.
    %\end{split}
%\end{equation}
%If $0 < \gamma \leq 1/2$, we will show that the right-hand side of \eqref{eq: limit double sum final} still converges. Note that by \eqref{eq: J_2 asymptotic},
%\begin{equation*}
     %\frac{J_{\gamma + 1}\big(2\sqrt{E(a^2m^2 + b^2n^2)}\big)}{(a^2m^2 + b^2n^2)^{(\gamma + 1)/2}} = C_{E,\gamma} \frac{\cos\big(2\sqrt{E(a^2m^2+b^2n^2)} - \frac{(2\gamma + 3)\pi}{4}\big)}{(a^2m^2 + b^2n^2)^{\gamma/2 + 3/4}} + O((a^2m^2+b^2n^2)^{-\gamma/2-1}).
%\end{equation*}
%Since $(a^2m^2+b^2n^2)^{-\gamma/2-1}$ is absolutely summable over $m,n \geq 1$, it is equivalent to show that
%\begin{equation*}
    %\lim_{N \to \infty} \sum_{1 \leq m,n \leq N} \frac{\cos\big(2\sqrt{E(a^2m^2+b^2n^2)} - \frac{(2\gamma + 3)\pi}{4}\big)}{(a^2m^2 + b^2n^2)^{\gamma/2 + 3/4}}
%\end{equation*}
%exists. Moreover, the constant phase shift is harmless, so it suffices to prove the following.

Putting everything together in \eqref{find asymptotic for}, we obtain the Riesz mean
\begin{equation} \label{eq: Riesz mean rectangle Bessel functions}
    \begin{split}
        R_\square^\gamma(E) &= \frac{ab}{4\pi(\gamma + 1)} E^{\gamma + 1} \pm \frac{a + b}{4\sqrt{\pi}}\frac{\Gamma(\gamma + 1)}{\Gamma(\gamma + 3/2)} E^{\gamma + 1/2} + \frac{1}{4} E^\gamma \\
        &+ \frac{\Gamma(\gamma+1)}{2\pi}E^{(\gamma+1)/2}\bigg(2ab \sum_{m,n \geq 1}^\square \frac{J_{\gamma + 1}\big(2\sqrt{E(a^2m^2 + b^2n^2)}\big)}{(a^2m^2 + b^2n^2)^{(\gamma + 1)/2}} \\
        &+ \frac{b}{a^\gamma}  \sum_{m = 1}^\infty \frac{J_{\gamma + 1}(2am\sqrt{E})}{m^{\gamma + 1}} + \frac{a}{b^\gamma} \sum_{n = 1}^\infty \frac{J_{\gamma+1}(2 b n \sqrt{E})}{n^{\gamma + 1}}\bigg) \\
        &\pm \frac{\Gamma(\gamma + 1)}{2\sqrt{\pi}}E^{\gamma/2 + 1/4}\bigg(\frac{1}{a^{\gamma-1/2}}  \sum_{m = 1}^\infty \frac{J_{\gamma + 1/2}(2am\sqrt{E})}{m^{\gamma+1/2}} + \frac{1}{b^{\gamma-1/2}} \sum_{n = 1}^\infty \frac{J_{\gamma + 1/2}(2bn\sqrt{E})}{n^{\gamma+1/2}}\bigg), 
    \end{split}
\end{equation}
for $\gamma > 0$. Substituting \eqref{eq: J_2 asymptotic} into \eqref{eq: Riesz mean rectangle Bessel functions} gives the expansion
\begin{equation} \label{eq: Riesz mean rectangle asymptotic}
    \begin{split}
        R_\square^\gamma(E) &= \frac{ab}{4\pi(\gamma+1)}E^{\gamma+1} \pm \frac{a+b}{4\sqrt{\pi}} \frac{\Gamma(\gamma+1)}{\Gamma(\gamma+3/2)}E^{\gamma+1/2} +\frac{1}{4}E^\gamma \\
        &+ \frac{\Gamma(\gamma + 1)}{\pi^{3/2}}E^{\gamma/2+1/4}\bigg(ab \sum_{m,n \geq 1}^\square \frac{\cos\big(2\sqrt{E(a^2m^2+b^2n^2)} -\frac{(2\gamma+3)\pi}{4}\big)}{(a^2m^2+b^2n^2)^{(2\gamma+3)/4}} \\
        &+ \frac{b}{2a^{\gamma+1/2}}\sum_{m=1}^\infty \frac{\cos\big(2am\sqrt{E} -\frac{(2\gamma+3)\pi}{4}\big)}{m^{\gamma+3/2}} + \frac{a}{2b^{\gamma+1/2}} \sum_{n=1}^\infty \frac{\cos\big(2bn\sqrt{E}-\frac{(2\gamma+3)\pi}{4}\big)}{n^{\gamma+3/2}}\bigg) \\
        &\pm \frac{\Gamma(\gamma+1)}{2\pi} E^{\gamma/2}\bigg(\frac{1}{a^\gamma} \sum_{m=1}^\infty
        \frac{\cos\big(2am\sqrt{E} - \frac{(\gamma+1)\pi}{2}\big)}{m^{\gamma+1}} + \frac{1}{b^\gamma} \sum_{n=1}^\infty \frac{\cos\big(2bn\sqrt{E} - \frac{(\gamma+1)\pi}{2}\big)
        }{n^{\gamma+1}}\bigg) \\
        &+ O(E^{\gamma/2 - 1/4}), \quad E \to \infty.
    \end{split}
\end{equation}
Apart from the double sum, which converges absolutely for $\gamma > 1/2$ and under rectangular summation for $\gamma > 0$, every sum in \eqref{eq: Riesz mean rectangle asymptotic} converges absolutely for $\gamma > 0$. In Appendix \ref{double sum bounds appendix}, we establish upper bounds for the double sum in \eqref{eq: Riesz mean rectangle asymptotic} and show that it is unbounded below for $0 < \gamma \leq 1/2$. By including further terms in the asymptotic expansion \eqref{eq: J_2 asymptotic} of $J_\nu$, one obtains the asymptotic expansion of $R_\square^\gamma(E)$ to arbitrary order. Since the arguments of $J_\nu$ in \eqref{eq: Riesz mean rectangle Bessel functions} are all proportional to $\sqrt{E}$, each successive term in the asymptotic expansion contributes an additional factor $E^{-1/2}$. Moreover, \eqref{eq: Riesz mean rectangle Bessel functions} involves both $J_{\gamma+1}$ and $J_{\gamma+1/2}$, whose leading contributions are of orders $E^{\gamma/2+1/4}$ and $E^{\gamma/2}$, respectively. The two resulting sequences of exponents therefore combine to give the powers $E^{(\gamma+1)/2-k/4}$ appearing in Theorem \ref{Riesz mean thm}. Each successive term in the asymptotic expansion of $J_\nu$ also gains an additional factor $(a^2m^2 + b^2n^2)^{-1/2}$ in the double sums, or an additional factor $m^{-1}$ or $n^{-1}$ in the single sums. Consequently, every series arising beyond the leading double series contribution converges absolutely for every $\gamma > 0$. Finally, since the full asymptotic expansion of $J_\nu$ consists of alternating cosine and sine terms, all coefficient functions are trigonometric series. This, together with Lemma \ref{shifted double sum bounds proposition}, proves Theorem \ref{Riesz mean thm} in the case of rectangles.

\begin{remark} It is interesting to note that in the case $a = b$, \eqref{eq: Riesz mean rectangle Bessel functions} is related to the Gauss circle problem. Indeed, letting $\gamma \to 0^+$, the Riesz mean becomes the eigenvalue counting function
\begin{equation*}
    R_\square^\gamma(E) \to R_\square^0(E) = \#\bigg\{(m,n) \in \NN^2 : m^2 + n^2 < \frac{a^2}{\pi^2}E\bigg\},
\end{equation*}
and the double sum in \eqref{eq: Riesz mean rectangle Bessel functions} formally approaches the series in Hardy's formula for the circle discrepancy \cite[Section 3]{ha15}.
\end{remark}

\section{Riesz mean asymptotics for integrable triangles} \label{triangles section}
In this section, we explain how the proof strategy of Section \ref{prelim section} yields the Riesz means for the integrable triangles. Recall that the Laplace eigenvalues for the equilateral triangle with side length $\ell$, the hemi-equilateral triangle with hypotenuse of length $\ell$, and the isosceles right triangle with legs of length $\ell$ with the Dirichlet boundary condition are
\begin{align*}
    \lambda^{D}_{m,n}(\triangledown) &= \frac{16 \pi^2}{9 \ell^2} (m^2 + mn + n^2), \qquad m, n > 0, \\
    \lambda_{m,n}^D(\heartsuit) &= \frac{16\pi^2}{9\ell^2}(m^2+mn+n^2), \qquad m > n > 0, \\
    \lambda^D_{m,n}(\diamondsuit) &= \frac{\pi^2}{\ell^2}(m^2+n^2), \qquad  m > n > 0,
\end{align*}
respectively, see \cite{maro25}. In each case, the eigenvalues with the Neumann boundary condition are obtained by replacing the strict inequalities involving $m$ and $n$ by $\geq$. In \cite{maro25}, the following identities for the corresponding heat traces are established:
\begin{equation} \label{heat traces triangles}
    \begin{split}
        H_\triangledown(t) &= \frac{\ell^2\sqrt{3}}{16\pi t} \pm \frac{3\ell}{8\sqrt{\pi t}} + \frac{1}{3} \pm \frac{3\ell}{4\sqrt{\pi t}}\sum_{m=1}^\infty e^{-9\ell^2m^2/(16t)} \\
        &+\frac{\ell^2\sqrt{3}}{8\pi t}\sum_{m=1}^\infty e^{-9\ell^2m^2/(4t)} + \frac{\ell^2 \sqrt{3}}{8 \pi t} \sum_{n=1}^\infty e^{-3\ell^2 n^2/(4t)} \\
        &+\frac{\ell^2\sqrt{3}}{4\pi t}\sum_{m,n\geq 1}e^{-3\ell^2(3m^2+n^2)/(4t)} + \frac{\ell^2\sqrt{3}}{4\pi t}\sum_{m,n\geq 1}e^{-3\ell^2(3(m-1/2)^2+(n-1/2)^2)/(4t)}, \\
        H_\heartsuit(t) &= \frac{\ell^2\sqrt{3}}{32\pi t} \pm \frac{\ell(3+\sqrt{3})}{16\sqrt{\pi t}} + \frac{5}{12} \pm \frac{\ell\sqrt{3}}{8\sqrt{\pi t}}\sum_{m=1}^\infty e^{-3\ell^2m^2/(16t)} \pm \frac{3\ell}{8\sqrt{\pi t}}\sum_{m=1}^\infty e^{-9\ell^2m^2/(16t)} \\
        &+\frac{\ell^2\sqrt{3}}{16\pi t}\sum_{n=1}^\infty e^{-3\ell^2n^2/(4t)} + \frac{\ell^2\sqrt{3}}{16\pi t}\sum_{m=1}^\infty e^{-9\ell^2m^2/(4t)} \\
        &+ \frac{\ell^2\sqrt{3}}{8\pi t}\sum_{m,n \geq 1} e^{-3\ell^2(3m^2+n^2)/(4t)} + \frac{\ell^2\sqrt{3}}{8\pi t}\sum_{m,n \geq 1} e^{-3\ell^2(3(m-1/2)^2+(n-1/2)^2)/(4t)}, \\
        H_\diamondsuit(t) &= \frac{\ell^2}{8\pi t} \pm \frac{\ell(2+\sqrt{2})}{8\sqrt{\pi t}} + \frac{3}{8} \pm \frac{\ell\sqrt{2}}{4\sqrt{\pi t}}\sum_{m=1}^\infty e^{-m^2\ell^2/(2t)} \pm \frac{\ell}{2\sqrt{\pi t}}\sum_{m=1}^\infty e^{-m^2\ell^2/t} \\ 
        &+ \frac{\ell^2}{2\pi t}\sum_{m=1}^\infty e^{-m^2\ell^2/t} + \frac{\ell^2}{2\pi t}\sum_{m,n \geq 1} e^{-\ell^2(m^2 + n^2)/t}, \quad t > 0.
    \end{split}
\end{equation}
For every term with $\pm$, we choose plus for the Neumann boundary condition and minus for the Dirichlet boundary condition. As in Section \ref{prelim section}, we analytically continue the heat traces to a complex argument $z$ with $\operatorname{Re}(z) > 0$, and replace $t$ by $\sigma + it, \,\, \sigma > 0, \,\, t \in \R$. As before, we multiply by $\varphi(t)$, integrate over $t \in \R$, and send $\sigma \to 0^+$ to obtain the Riesz means
\begin{equation*}
    \begin{split}
        R_\triangledown^\gamma(E) &= \frac{\ell^2\sqrt{3}}{16\pi(\gamma + 1)}E^{\gamma + 1} \pm \frac{3\ell}{8\sqrt{\pi}}\frac{\Gamma(\gamma+1)}{\Gamma(\gamma + 3/2)}E^{\gamma +1/2} + \frac{1}{3}E^\gamma + \lim_{\sigma \to 0^+} \int_{-\infty}^\infty \varphi(t) \mathcal R_\triangledown(\sigma + it)dt, \\
        R_\heartsuit^\gamma(E) &= \frac{\ell^2\sqrt{3}}{32\pi(\gamma + 1)}E^{\gamma + 1} \pm \frac{\ell(3 + \sqrt{3})}{16\sqrt{\pi}}\frac{\Gamma(\gamma+1)}{\Gamma(\gamma + 3/2)}E^{\gamma +1/2} + \frac{5}{12}E^\gamma + \lim_{\sigma \to 0^+} \int_{-\infty}^\infty \varphi(t) \mathcal R_\heartsuit(\sigma + it)dt, \\
        R_\diamondsuit^\gamma(E) &= \frac{\ell^2}{8\pi(\gamma + 1)}E^{\gamma + 1} \pm \frac{\ell(2 + \sqrt{2})}{8\sqrt{\pi}}\frac{\Gamma(\gamma+1)}{\Gamma(\gamma + 3/2)}E^{\gamma + 1/2} + \frac{3}{8}E^\gamma + \lim_{\sigma \to 0^+} \int_{-\infty}^\infty \varphi(t) \mathcal R_\diamondsuit(\sigma + it)dt,
    \end{split}
\end{equation*}
where
\begin{equation*}
    \begin{split}
        \mathcal{R}_\triangledown(\sigma + it) &= \pm \frac{3\ell}{4\sqrt{\pi (\sigma + it)}}\sum_{m=1}^\infty e^{-9\ell^2m^2/(16(\sigma + it))} + \frac{\ell^2\sqrt{3}}{8\pi (\sigma + it)} \sum_{m=1}^\infty e^{-9\ell^2m^2/(4(\sigma + it))} \\
        &+ \frac{\ell^2 \sqrt{3}}{8 \pi (\sigma + it)} \sum_{n=1}^\infty e^{-3\ell^2 n^2/(4(\sigma + it))} + \frac{\ell^2\sqrt{3}}{4\pi (\sigma + it)} \sum_{m,n\geq 1} e^{-3\ell^2(3m^2+n^2)/(4(\sigma + it))} \\
        &+ \frac{\ell^2\sqrt{3}}{4\pi (\sigma + it)}\sum_{m,n\geq 1} e^{-3\ell^2(3(m-1/2)^2+(n-1/2)^2)/(4(\sigma + it))}, \\
        \mathcal{R}_\heartsuit(\sigma + it) &= \pm \frac{\ell\sqrt{3}}{8\sqrt{\pi (\sigma + it)}}\sum_{m=1}^\infty e^{-3\ell^2m^2/(16(\sigma + it))} \pm \frac{3\ell}{8\sqrt{\pi (\sigma + it)}}\sum_{m=1}^\infty e^{-9\ell^2m^2/(16(\sigma + it))} \\
        &+\frac{\ell^2\sqrt{3}}{16\pi (\sigma + it)}\sum_{n=1}^\infty e^{-3\ell^2n^2/(4(\sigma + it))} + \frac{\ell^2\sqrt{3}}{16\pi (\sigma + it)}\sum_{m=1}^\infty e^{-9\ell^2m^2/(4(\sigma + it))} \\
        &+ \frac{\ell^2\sqrt{3}}{8\pi (\sigma + it)}\sum_{m,n \geq 1} e^{-3\ell^2(3m^2+n^2)/(4(\sigma + it))} + \frac{\ell^2\sqrt{3}}{8\pi (\sigma + it)}\sum_{m,n \geq 1} e^{-3\ell^2(3(m-1/2)^2+(n-1/2)^2)/(4(\sigma + it))}, \\
        \mathcal{R}_\diamondsuit(\sigma + it) &= \pm \frac{\ell\sqrt{2}}{4\sqrt{\pi (\sigma + it)}}\sum_{m=1}^\infty e^{-m^2\ell^2/(2(\sigma + it))} \pm \frac{\ell}{2\sqrt{\pi (\sigma + it)}}\sum_{m=1}^\infty e^{-m^2\ell^2/(\sigma + it)} \\ 
        &+ \frac{\ell^2}{2\pi (\sigma + it)}\sum_{m=1}^\infty e^{-m^2\ell^2/(\sigma + it)} + \frac{\ell^2}{2\pi (\sigma + it)}\sum_{m,n \geq 1} e^{-\ell^2(m^2 + n^2)/(\sigma + it)}.
    \end{split}
\end{equation*}
Note that in addition to the sums appearing in Section \ref{prelim section}, here we also have sums involving half-integers. The same proof as in Lemma \ref{lemma: sharp exponential bound on partial sums} shows that
\begin{equation*} \label{eq: sharp exponential bound partial sums}
    \sum_{n = 1}^N e^{-(n-1/2)^2/(\sigma + it)}  = O_\sigma((1 + |t|)^{1/2}),
\end{equation*}
%The proof is the same as in
%Lemma~\ref{lemma: sharp exponential bound on partial sums}.
%Writing $(\sigma+it)^{-1}=\alpha-i\beta$, we note that the
%shift leaves the second derivative of the phase unchanged:
%$\bigl(\beta(x-\tfrac12)^2\bigr)''=2\beta$. Thus, for $|t|\geq1$,
%van der Corput gives
%\begin{equation*}
    %\left|\sum_{n=1}^M e^{i\beta(n-\frac12)^2}\right|
    %\lesssim M|\beta|^{1/2}+|\beta|^{-1/2}.
%\end{equation*}
%The weights $w_n=e^{-\alpha(n-\frac12)^2}$
%satisfy $\sum_{n\geq1}w_n\lesssim\alpha^{-1/2} \lesssim \alpha^{-1/2}$. This is shown by using $n-\tfrac12\geq n/2$ and then estimating the sum by a Gaussian integral. Summation by parts therefore bounds $\sum_{n=1}^N e^{-\alpha(n - \frac12}^2} e^{i\beta(n - \frac12)^2}$, uniformly in $N$, by
%\begin{equation*}
    %C\left(|\beta|^{1/2}\sum_{n=1}^N w_n
           %+|\beta|^{-1/2}w_1\right)
    %\lesssim \sqrt{\frac{|\beta|}{\alpha}}+|\beta|^{-1/2}
    %\lesssim_\sigma |t|^{1/2}.
%\end{equation*}
%For $|t|\leq1$, the Gaussian weight is summable uniformly in $t$,
%since $\alpha\geq\sigma/(\sigma^2+1)$.
uniformly in $N \in \NN$ and $t \in \R$. Consequently, all interchanges between integrals and infinite sums are justified exactly as in Section \ref{prelim section}.
Thus, we can compute
\begin{equation*}
    \lim_{\sigma \to 0^+} \int_{-\infty}^\infty \varphi(t) \mathcal R_\triangledown(\sigma + it)dt, \quad \lim_{\sigma \to 0^+} \int_{-\infty}^\infty \varphi(t) \mathcal R_\heartsuit(\sigma + it)dt, \quad \lim_{\sigma \to 0^+} \int_{-\infty}^\infty \varphi(t) \mathcal R_\diamondsuit(\sigma + it)dt
\end{equation*}
exactly as in %\eqref{eq: DCT -1/2}, \eqref{eq single sum -1/2}, \eqref{eq: DCT -1}, \eqref{eq single sum -1}, \eqref{eq: DCT double}, we obtain
\eqref{eq: DCT -1/2}-\eqref{eq: limit double sum final} to obtain
\begin{equation} \label{eq: Riesz mean equilateral Bessel functions}
    \begin{split}
        R_\triangledown^\gamma(E) &= \frac{\ell^2\sqrt{3}}{16\pi(\gamma+1)}E^{\gamma+1} \pm \frac{3\ell}{8\sqrt{\pi}} \frac{\Gamma(\gamma+1)}{\Gamma(\gamma+3/2)} E^{\gamma+1/2}
        +\frac13E^\gamma \\
        &+ \frac{2^{\gamma-1}\Gamma(\gamma+1)}{3^{\gamma/2}\pi\,\ell^{\gamma-1}}E^{(\gamma+1)/2} \bigg(\sum_{m,n \geq 1}^\square \frac{J_{\gamma+1}\big(\ell\sqrt{3E(3m^2+n^2)}\big)}{(3m^2+n^2)^{(\gamma+1)/2}} \\
        &+ \sum_{m,n \geq 1}^\square \frac{J_{\gamma+1}\big(\ell\sqrt{3E(3(m-1/2)^2+(n-1/2)^2)}\big)}{\big(3(m-1/2)^2+(n-1/2)^2\big)^{(\gamma+1)/2}} \\
        &+ \frac{1}{2 \cdot 3^{(\gamma+1)/2}} \sum_{m=1}^\infty \frac{J_{\gamma+1}(3\ell m\sqrt E)}{m^{\gamma+1}} + \frac{1}{2} \sum_{n=1}^\infty \frac{J_{\gamma+1}(\ell\sqrt3\,n\sqrt E)}{n^{\gamma+1}}\bigg) \\
        &\pm \frac{4^{\gamma-1/2}\Gamma(\gamma+1)}{3^{\gamma-1/2}\sqrt{\pi}\,\ell^{\gamma-1/2}}E^{\gamma/2+1/4} \sum_{m=1}^\infty \frac{J_{\gamma+1/2}(3\ell m\sqrt E/2)}{m^{\gamma+1/2}},
    \end{split}
\end{equation}
\begin{equation} \label{eq: Riesz mean hemiequilateral Bessel functions}
    \begin{split}
        R_\heartsuit^\gamma(E) &= \frac{\ell^2\sqrt3}{32\pi(\gamma+1)}E^{\gamma+1} \pm \frac{\ell(3+\sqrt3)}{16\sqrt\pi} \frac{\Gamma(\gamma+1)}{\Gamma(\gamma+3/2)}E^{\gamma+1/2} +\frac{5}{12}E^\gamma \\
        &+ \frac{2^{\gamma-2}\Gamma(\gamma+1)}{3^{\gamma/2}\pi\,\ell^{\gamma-1}} E^{(\gamma+1)/2}\bigg(
        \sum_{m,n \geq 1}^\square \frac{J_{\gamma+1}\big(\ell\sqrt{3E(3m^2+n^2)}\big)}{(3m^2+n^2)^{(\gamma+1)/2}} \\
        &+ \sum_{m,n \geq 1}^\square \frac{J_{\gamma+1}\big(\ell\sqrt{3E(3(m-1/2)^2+(n-1/2)^2)}\big)}{\big(3(m-1/2)^2+(n-1/2)^2\big)^{(\gamma+1)/2}}\bigg) \\
        &+ \frac{2^{\gamma-3}\Gamma(\gamma+1)}{3^{\gamma/2}\pi\,\ell^{\gamma-1}} E^{(\gamma+1)/2}\bigg(\frac{1}{3^{(\gamma+1)/2}} \sum_{m=1}^\infty \frac{J_{\gamma+1}(3\ell m\sqrt E)}{m^{\gamma+1}} + \sum_{n=1}^\infty \frac{J_{\gamma+1}(\ell\sqrt{3}\, n\sqrt E)}{n^{\gamma+1}}\bigg) \\
        &\pm \frac{2^{2\gamma-2}\Gamma(\gamma+1)}{3^{\gamma/2-1/4}\ell^{\gamma-1/2}\sqrt{\pi}} E^{\gamma/2+1/4}
        \bigg(\sum_{m=1}^\infty \frac{J_{\gamma+1/2}\big(\ell\sqrt{3}\,m\sqrt E/2\big)}{m^{\gamma+1/2}} \\
        &+ \frac{1}{3^{\gamma/2-1/4}} \sum_{m=1}^\infty \frac{J_{\gamma+1/2}\big(3\ell m\sqrt E/2\big)}{m^{\gamma+1/2}}\bigg),
    \end{split}
\end{equation}
\begin{equation} \label{eq: Riesz mean isosceles Bessel functions}
    \begin{split}
        R_\diamondsuit^\gamma(E) &= \frac{\ell^2}{8\pi(\gamma+1)}E^{\gamma+1}
        \pm \frac{\ell(2+\sqrt2)}{8\sqrt\pi}\frac{\Gamma(\gamma+1)}{\Gamma(\gamma+3/2)}E^{\gamma+1/2} + \frac{3}{8}E^\gamma \\
        &+ \frac{\Gamma(\gamma+1)}{2\pi\,\ell^{\gamma-1}}E^{(\gamma+1)/2} \bigg(\sum_{m,n \geq 1}^\square
        \frac{J_{\gamma+1}\big(2\ell\sqrt{E(m^2+n^2)}\big)}
        {(m^2+n^2)^{(\gamma+1)/2}} + \sum_{m=1}^\infty
        \frac{J_{\gamma+1}(2\ell m\sqrt E)}{m^{\gamma+1}}\bigg) \\
        &\pm \frac{\Gamma(\gamma+1)}{\ell^{\gamma-1/2}\sqrt{\pi}}E^{\gamma/2+1/4}\bigg(2^{\gamma/2-5/4} \sum_{m=1}^\infty \frac{J_{\gamma+1/2}(\ell \sqrt{2}\, m\sqrt E)}{m^{\gamma+1/2}} + \frac{1}{2} \sum_{m=1}^\infty \frac{J_{\gamma+1/2}(2\ell m\sqrt E)}{m^{\gamma+1/2}}\bigg),
    \end{split}
\end{equation}
for $\gamma > 0$. Using \eqref{eq: J_2 asymptotic}, these become
\begin{equation} \label{eq: Riesz mean equilateral asymptotic}
    \begin{split}
        R_\triangledown^\gamma(E) &= \frac{\ell^2\sqrt{3}}{16\pi(\gamma+1)}E^{\gamma+1} \pm \frac{3\ell}{8\sqrt{\pi}} \frac{\Gamma(\gamma+1)}{\Gamma(\gamma+3/2)} E^{\gamma+1/2} +\frac{1}{3}E^\gamma \\
        &+ \frac{2^{\gamma-3/2}\Gamma(\gamma+1)}{3^{\gamma/2+1/4}\pi^{3/2}\ell^{\gamma-1/2}} E^{\gamma/2+1/4}\bigg(2 \sum_{m,n \geq 1}^\square \frac{\cos\big(\ell\sqrt{3E(3m^2+n^2)}-\frac{(2\gamma+3)\pi}{4}\big)}{(3m^2+n^2)^{(2\gamma+3)/4}} \\
        &+ 2 \sum_{m,n \geq 1}^\square \frac{\cos\big(\ell\sqrt{3E(3(m-1/2)^2+(n-1/2)^2)}-\frac{(2\gamma+3)\pi}{4}\big)}{\big(3(m-1/2)^2+(n-1/2)^2\big)^{(2\gamma+3)/4}} \\
        &+ \frac{1}{3^{\gamma/2+3/4}} \sum_{m=1}^\infty \frac{\cos\big(3\ell m\sqrt E-\frac{(2\gamma+3)\pi}{4}\big)}{m^{\gamma+3/2}} + \sum_{n=1}^\infty \frac{\cos\big(\ell\sqrt3\,n\sqrt E-\frac{(2\gamma+3)\pi}{4}\big)}{n^{\gamma+3/2}}\bigg) \\
        &\pm \frac{4^{\gamma}\Gamma(\gamma+1)}{3^\gamma\pi\ell^\gamma} E^{\gamma/2} \sum_{m=1}^\infty \frac{\cos\big(3\ell m\sqrt{E}/2-\frac{(\gamma+1)\pi}{2}\big)}{m^{\gamma+1}} + O(E^{\gamma/2-1/4}),
    \end{split}
\end{equation}
\begin{equation} \label{eq: Riesz mean hemiequilateral asymptotic}
    \begin{split}
        R_\heartsuit^\gamma(E) &= \frac{\ell^2\sqrt3}{32\pi(\gamma+1)}E^{\gamma+1} \pm \frac{\ell(3+\sqrt3)}{16\sqrt\pi} \frac{\Gamma(\gamma+1)}{\Gamma(\gamma+3/2)} E^{\gamma+1/2} + \frac{5}{12}E^\gamma \\
        &+ \frac{2^{\gamma-5/2}\Gamma(\gamma+1)}{3^{\gamma/2+1/4}\pi^{3/2}\ell^{\gamma-1/2}}E^{\gamma/2+1/4} \bigg(2 \sum_{m,n \geq 1}^\square \frac{\cos\big(\ell\sqrt{3E(3m^2+n^2)} - \frac{(2\gamma+3)\pi}{4}\big)}{(3m^2+n^2)^{(2\gamma+3)/4}} \\
        &+ 2\sum_{m,n \geq 1}^\square \frac{\cos\big(\ell\sqrt{3E(3(m-1/2)^2+(n-1/2)^2)} -\frac{(2\gamma+3)\pi}{4}\big)}{\big(3(m-1/2)^2+(n-1/2)^2\big)^{(2\gamma+3)/4}} \\
        &+ \frac{1}{3^{\gamma/2+3/4}} \sum_{m=1}^\infty \frac{\cos\big(3\ell m\sqrt E-\frac{(2\gamma+3)\pi}{4}\big)}{m^{\gamma+3/2}} + \sum_{n=1}^\infty \frac{\cos\big(\ell\sqrt3\,n\sqrt E-\frac{(2\gamma+3)\pi}{4}\big)}{n^{\gamma+3/2}}\bigg) \\
        &\pm \frac{2^{2\gamma-1}\Gamma(\gamma+1)}{3^{\gamma/2}\pi\ell^\gamma}E^{\gamma/2}\bigg( \sum_{m=1}^\infty \frac{\cos\big(\ell\sqrt3\,m\sqrt E/2-\frac{(\gamma+1)\pi}{2}\big)}{m^{\gamma+1}} \\
        &+ \frac{1}{3^{\gamma/2}} \sum_{m=1}^\infty \frac{\cos\big(3\ell m\sqrt E/2-\frac{(\gamma+1)\pi}{2}\big)}{m^{\gamma+1}}\bigg) + O(E^{\gamma/2-1/4}),
    \end{split}
\end{equation}
\begin{equation} \label{eq: Riesz mean isosceles asymptotic}
    \begin{split}
        R_\diamondsuit^\gamma(E) &= \frac{\ell^2}{8\pi(\gamma+1)}E^{\gamma+1} \pm \frac{\ell(2+\sqrt2)}{8\sqrt\pi} \frac{\Gamma(\gamma+1)}{\Gamma(\gamma+3/2)}E^{\gamma+1/2} + \frac{3}{8}E^\gamma \\
        &+ \frac{\Gamma(\gamma+1)}{2\pi^{3/2}\ell^{\gamma-1/2}}E^{\gamma/2+1/4} \bigg(\sum_{m,n \geq 1}^\square \frac{
        \cos\big(2\ell\sqrt{E(m^2+n^2)}-\frac{(2\gamma+3)\pi}{4}\big)}{(m^2+n^2)^{(2\gamma+3)/4}} \\
        &+ \sum_{m=1}^\infty \frac{
        \cos\big(2\ell m\sqrt E-\frac{(2\gamma+3)\pi}{4}\big)}{m^{\gamma+3/2}}\bigg) \\
        &\pm \frac{\Gamma(\gamma+1)}{2\pi\ell^\gamma}E^{\gamma/2}\bigg(2^{\gamma/2} \sum_{m=1}^\infty \frac{
        \cos\big(\ell\sqrt2\,m\sqrt E-\frac{(\gamma+1)\pi}{2}\big)}{m^{\gamma+1}} + \sum_{m=1}^\infty \frac{\cos\big(2\ell m\sqrt E-\frac{(\gamma+1)\pi}{2}\big)}{m^{\gamma+1}}\bigg) \\
        &+ O(E^{\gamma/2-1/4}).
    \end{split}
\end{equation}
Apart from the double sums, which converge absolutely for $\gamma > 1/2$ and under rectangular summation for $\gamma > 0$, every sum in \eqref{eq: Riesz mean equilateral asymptotic}-\eqref{eq: Riesz mean isosceles asymptotic} converges absolutely for every $\gamma > 0$. By including further terms in the asymptotic expansion \eqref{eq: J_2 asymptotic} of $J_\nu$, one obtains the expansions of the Riesz means to arbitrary order. Since the arguments of the Bessel functions in \eqref{eq: Riesz mean equilateral Bessel functions}-\eqref{eq: Riesz mean isosceles Bessel functions} are all proportional to $\sqrt{E}$, each successive term contributes an additional factor $E^{-1/2}$. Moreover, these formulas involve both $J_{\gamma + 1}$ and $J_{\gamma + 1/2}$, whose leading contributions are of orders $E^{\gamma/2 + 1/4}$ and $E^{\gamma/2}$, respectively. The two resulting sequences of exponents therefore combine to give the powers $E^{(\gamma+1)/2 - k/4}$ appearing in Theorem \ref{Riesz mean thm}. Each successive term in the Bessel asymptotic expansion also gains an additional factor of the inverse radial variable in the corresponding series. Consequently, every series arising beyond the leading double-series contributions converges absolutely for every $\gamma > 0$. Finally, since the full asymptotic expansion of $J_\nu$ consists of alternating cosine and sine terms, all coefficient functions are trigonometric series. This, together with Appendix \ref{double sum bounds appendix}, completes the proof of Theorem \ref{Riesz mean thm}.
%\newpage
\appendix

\section{Uniform integral estimates} \label{estimate appendix new}
We establish certain integral estimates used in Section \ref{prelim section}. 
Throughout, $E > 0$ and $\gamma > 0$ are fixed, and all implicit constants are uniform in $0 \leq \sigma \leq 1$.

To justify \eqref{eq: integral cosine estimate}, 
it suffices to show
\begin{equation}\label{cosine decay}
    I_\sigma(k) \defeq \int_0^{\sqrt E} (E - x^2)^\gamma e^{-\sigma x^2} \cos(kx) dx = O(k^{-\gamma - 1} ), \quad k \to \infty.
\end{equation}
Write
\begin{equation*}
    e^{-\sigma x^2} = e^{-\sigma E} e^{\sigma(E-x^2)} = e^{-\sigma E} \sum_{j = 0}^\infty \frac{\sigma^j(E-x^2)^j}{j!}.
\end{equation*}
Since the series converges uniformly for $0 \leq x \leq \sqrt{E}$ and $0 \leq \sigma \leq 1$, we may integrate termwise:
\begin{equation*}
    \begin{split}
        I_\sigma(k) &= e^{-\sigma E} \sum_{j = 0}^\infty \frac{\sigma^j}{j!} \int_0^{\sqrt{E}} (E - x^2)^{\gamma + j} \cos(kx) dx \\
        &= e^{-\sigma E} \sum_{j = 0}^\infty \frac{\sigma^j}{j!} \sqrt{\pi} \,2^{\gamma + j - 1/2} \, \Gamma(\gamma + j + 1) E^{\gamma/2 + j/2 + 1/4} k^{-\gamma-j-1/2} J_{\gamma + j + 1/2}(k\sqrt{E}).
    \end{split}
\end{equation*}
The last equality follows from \cite[\S 10.9(i), Equation 10.9.4]{dlmf}. By \eqref{eq: J_2 asymptotic}, the $j = 0$ term becomes
\begin{equation} \label{eq j = 0 1}
    e^{-\sigma E}  \, 2^{\gamma} \,\Gamma(\gamma + 1) E^{\gamma/2} k^{-\gamma-1} \cos\bigg(k\sqrt{E} - \frac{(\gamma + 1)\pi}{2}\bigg) + O(k^{-\gamma-2}).
\end{equation}
For $j \geq 1$, we use $|J_\nu(x)| \leq 1$ for $x \geq 0$ and $\nu \geq 0$ \cite[\S 10.14, equation 10.14.1]{dlmf}, together with
\begin{equation*}
    \frac{\Gamma(\gamma + j + 1)}{j!} = O(j^\gamma), \quad j \to \infty.
\end{equation*}
For $k \geq 4\sqrt{E}$, this gives
\begin{equation} \label{eq j => 1 1}
    \begin{split}
        \bigg|e^{-\sigma E} \sum_{j = 1}^\infty \frac{\sigma^j}{j!} \sqrt{\pi} \,2^{\gamma + j - 1/2} \, \Gamma(\gamma + j + 1) E^{\gamma/2 + j/2 + 1/4} k^{-\gamma-j-1/2} J_{\gamma + j + 1/2}(k\sqrt{E})\bigg| \\
        \lesssim k^{-\gamma-1/2} \sum_{j = 1}^\infty \frac{\Gamma(\gamma + j + 1)}{j!} \bigg(\frac{2\sqrt{E}}{k}\bigg)^j \lesssim k^{-\gamma - 3/2} \sum_{j = 1}^\infty j^\gamma 2^{-j} \lesssim k^{-\gamma - 3/2}.
    \end{split}
\end{equation}
Combining \eqref{eq j = 0 1} and \eqref{eq j => 1 1} gives
\begin{equation*}
    I_\sigma(k) = e^{-\sigma E}  \, 2^{\gamma} \,\Gamma(\gamma + 1) E^{\gamma/2} k^{-\gamma-1} \cos\bigg(k\sqrt{E} - \frac{(\gamma + 1)\pi}{2}\bigg) + O(k^{-\gamma-3/2}),
\end{equation*}
which proves \eqref{cosine decay}.

We next address the Bessel integrals in \eqref{eq single sum -1},
%\eqref{eq2 single sum -1}, and %\eqref{eq double sum}
\eqref{eq: integral J_0 estimate}, and \eqref{eq: integral h_sigma estimate}.
%\eqref{eq: limit double sum final}. 
We claim that
\begin{equation} \label{bessel bound}
    \begin{split}
    J_\sigma(k) &\defeq \int_0^{\sqrt{E}} x (E - x^2)^\gamma e^{-\sigma x^2} J_0(kx) dx = O(k^{-\gamma-3/2}), \quad k \to \infty, \\
    \Tilde{J}_\sigma(k) &\defeq \int_0^{\sqrt{E}} x (E - x^2)^\gamma h_\sigma(E - x^2) J_0(kx) dx = O(k^{-\gamma-3/2}), \quad k \to \infty,
    \end{split}
\end{equation}
where $h_\sigma(y) = \int_0^1 e^{s \sigma y} ds$. The proofs are entirely analogous to \eqref{cosine decay}, and follow from the calculations
\begin{equation} \label{eq: J_sigma first term}
    \begin{split}
        J_\sigma(k) &= e^{-\sigma E}\sum_{j = 0}^\infty \frac{\sigma^j}{j!} \int_0^{\sqrt{E}} x(E-x^2)^{\gamma + j} J_0(kx) dx \\
        &= e^{-\sigma E}\sum_{j = 0}^\infty \frac{\sigma^j}{j!} 2^{\gamma  + j} \Gamma(\gamma + j + 1) E^{(\gamma + j + 1)/2} k^{-\gamma - j - 1} J_{\gamma + j + 1}(k\sqrt{E}) \\
        &= e^{-\sigma E} 2^\gamma \sqrt{\frac{2}{\pi}} \Gamma(\gamma + 1)E^{\gamma/2 + 1/4} k^{-\gamma -3/2}\cos\bigg(k\sqrt{E} - \frac{(2\gamma  + 3)\pi}{4}\bigg) + O(k^{-\gamma - 2}), 
    \end{split}
\end{equation}
\begin{equation} \label{eq: J_sigma tilde first term}
    \begin{split}
        \Tilde{J}_\sigma(k) &= \sum_{j = 0}^\infty \frac{\sigma^j}{(j+1)!} \int_0^{\sqrt{E}} x(E-x^2)^{\gamma + j} J_0(kx) dx \\
        &= \sum_{j = 0}^\infty \frac{\sigma^j}{(j+1)!} 2^{\gamma + j} \Gamma(\gamma + j + 1) E^{(\gamma + j + 1)/2} k^{-\gamma-j-1} J_{\gamma+j+1}(k\sqrt{E}) \\
        &= 2^\gamma \sqrt{\frac{2}{\pi}} \Gamma(\gamma + 1) E^{\gamma/2 + 1/4} k^{-\gamma-3/2}\cos\bigg(k\sqrt{E} - \frac{(2\gamma  + 3)\pi}{4}\bigg) + O(k^{-\gamma-2}),
    \end{split}
\end{equation}
where we again used \cite{we}.

\section{Bounds for double sums} \label{double sum bounds appendix}
Fix $a, b > 0$ and $0 \le s \le 1/2$. Put
\begin{equation*}
    r_{m,n} \defeq \bigl(a^2(m-s)^2 + b^2(n-s)^2\bigr)^{1/2}.
\end{equation*}
For $\gamma > 0$, $t \geq 1$, and $N \in \NN$, define
\begin{equation*}
    S_{\gamma,s,N}(t) \defeq \sum_{1 \leq m,n \leq N} \frac{\cos(2tr_{m,n} - \theta_\gamma)} {r_{m,n}^{\gamma + 3/2}}, \quad \theta_\gamma \defeq \frac{(2\gamma + 3)\pi}{4},
\end{equation*}
and set
\begin{equation} \label{eq: S gamma s}
    S_{\gamma,s}(t) \defeq \lim_{N\to\infty}S_{\gamma,s,N}(t).
\end{equation}
We now establish upper bounds and show that $S_{\gamma,s}(t)$ is unbounded in $t$ for $0 < \gamma \leq 1/2$.
\begin{lemma} \label{shifted double sum bounds proposition}
    The limit $S_{\gamma,s}(t)$ exists for every $\gamma > 0$ and $t \geq 1$. %\textcolor{red}{Did we already show this limit exists in an earlier section?} 
    Moreover, uniformly in $N \in \NN$,
    \begin{equation} \label{shifted double sum upper}
        |S_{\gamma,s,N}(t)| \leq C
        \begin{cases}
            t^{1/2-\gamma}, & 0 < \gamma < \frac{1}{2} \\
            1 + \log t, & \gamma = \frac{1}{2} \\
            1, & \gamma > \frac{1}{2}
        \end{cases},
        \qquad t \geq 1.
    \end{equation}
    The partial sums converge uniformly on compact subsets of $[1, \infty)$.
    
    If $0 < \gamma \leq \frac{1}{2}$, there exist $c > 0$ and an increasing sequence $E_j \to \infty$ such that
    \begin{equation} \label{shifted double sum lower}
        -S_{\gamma,s}(\sqrt{E_j}) \geq c
        \begin{cases}
            (\log E_j)^{1/4-\gamma/2}, & 0 < \gamma < \frac{1}{2}, \\
            \log\log E_j, & \gamma = \frac{1}{2}.
        \end{cases}
    \end{equation}
\end{lemma}
\begin{proof}
%For notational brevity 
Put $p \defeq \gamma + \frac{3}{2}$. We first prove the upper bound \eqref{shifted double sum upper} by decomposing the rectangular sums into square shells. Set
\begin{equation*}
    T_N(t) \defeq \sum_{1 \leq m,n \leq N} \frac{e^{2itr_{m,n}}}{r_{m,n}^{p}}.
\end{equation*}
Then
\begin{equation*}
    T_N(t) = \sum_{k=1}^N H_k(t), \quad H_k(t) \defeq \sum_{m=1}^k \frac{e^{2itr_{m,k}}}{r_{m,k}^{p}} + \sum_{n=1}^{k-1} \frac{e^{2itr_{k,n}}}{r_{k,n}^{p}}.
\end{equation*}
Consider the first sum in $H_k(t)$. For $1 \leq M \leq k$, define
\begin{equation*}
    A_k(M) \defeq \sum_{m=1}^M e^{2itr_{m,k}}.
\end{equation*}
Clearly $|A_k(M)| \leq k$. On the other hand, for
\begin{equation*}
    f_k(x) \defeq 2t(a^2(x-s)^2+b^2(k-s)^2)^{1/2},
\end{equation*}
we have
\begin{equation*}
    f_k''(x) = \frac{2ta^2b^2(k-s)^2}{(a^2(x-s)^2+b^2(k-s)^2)^{3/2}} \asymp \frac{t}{k}, \quad 1 \leq x \leq k.
\end{equation*}
The van der Corput estimate \cite[Theorem 2.2]{grko91} therefore gives
\begin{equation*}
    |A_k(M)| \lesssim \bigg(k\sqrt{\frac{t}{k}} + \sqrt{\frac{k}{t}}\bigg) \lesssim \sqrt{tk}, \quad t \geq 1,
\end{equation*}
uniformly in $1 \leq M \leq k$. Combining the two estimates for $|A_k(M)|$ shows
\begin{equation*}
    |A_k(M)| \lesssim \min(k, \sqrt{tk}).
\end{equation*}

The weight
\begin{equation*}
    w_k(x) \defeq (a^2(x-s)^2+b^2(k-s)^2)^{-p/2}
\end{equation*}
is positive and decreasing in $x$, and $w_k(x)\lesssim k^{-p}$. Summation by parts gives
\begin{equation*}
    \sum_{m=1}^k\frac{e^{2itr_{m,k}}}{r_{m,k}^{p}} = A_k(k)w_k(k) + \sum_{m=1}^{k-1} A_k(m)(w_k(m)-w_k(m+1)).
\end{equation*}
Since
\begin{equation*}
    w_k(k) + \sum_{m=1}^{k-1} (w_k(m)-w_k(m+1)) = w_k(1) \lesssim k^{-p},
\end{equation*}
it follows that
\begin{equation*}
    \bigg|\sum_{m=1}^k \frac{e^{2itr_{m,k}}}{r_{m,k}^{p}}\bigg| \lesssim \min(k^{1-p},t^{1/2}k^{1/2-p}).
\end{equation*}
Interchanging $a$ and $b$ gives the same estimate for the second sum in $H_k(t)$. Thus,
\begin{equation}\label{eq:shifted-shell-bound}
    |H_k(t)| \lesssim \min(k^{-1/2-\gamma}, t^{1/2}k^{-1-\gamma}).
\end{equation}
In particular, $\sum_{k\geq1}|H_k(t)| < \infty$ for each fixed $t \geq 1$. Since 
\begin{equation*}
    S_{\gamma,s,N}(t) = \operatorname{Re}\bigl(e^{-i\theta_\gamma}T_N(t)\bigr),
\end{equation*}
this proves convergence of %the rectangular sums 
$S_{\gamma, s, N}(t)$ as $N \to \infty$. Furthermore, \eqref{eq:shifted-shell-bound} shows
\begin{equation*}
    \sum_{k=1}^{\infty} |H_k(t)| \lesssim \sum_{k \leq t} k^{-1/2-\gamma} + t^{1/2}\sum_{k > t} k^{-1-\gamma} \lesssim
    \begin{cases}
        t^{1/2-\gamma}, & 0 < \gamma < \frac{1}{2}, \\
        1 + \log t, & \gamma = \frac{1}{2}, \\
        1, & \gamma > \frac{1}{2}.
    \end{cases}
\end{equation*}
This confirms \eqref{shifted double sum upper}. Finally, for $0 < \alpha < \beta$, $\alpha \le t \le \beta$ and $N > M$,
\begin{equation} \label{eq: partial sums converge uniformly}
    |T_N(t) - T_M(t)| \leq \sum_{k = M+1}^N |H_k(t)| \lesssim \beta^{1/2}\sum_{k > M} k^{-1-\gamma} \lesssim \beta^{1/2}M^{-\gamma}.
\end{equation}
This establishes uniform convergence on compact subsets of $[1,\infty)$.

We now assume $0 < \gamma \leq \frac{1}{2}$ and prove the lower
bound \eqref{shifted double sum lower}. The strategy is as
follows. At scale $R^{-1}$, we smooth $S_{\gamma,s}$ with a
nonnegative kernel of integral one whose Fourier transform is
also nonnegative. The resulting series has nonnegative
coefficients, and its value at $t=0$ tends to $-\infty$ as
$R\to\infty$. For each large $R$, simultaneous approximation
gives a point $t_R\geq R$ at which the phases of the dominant
terms nearly return to their values at zero, making the average
comparably negative. Since the kernel is nonnegative, the
original sum must be at least as negative at some nearby point.

We begin by constructing the averaging kernel. Choose a nonnegative function $\psi\in C_0^\infty(\R)$ with $\supp \psi \subseteq [-1/2, 1/2]$ and
\begin{equation*}
    \int_{\R} \psi(u) du = 1.
\end{equation*}
Define
\begin{equation*}
    \rho(u) \defeq \int_{\R} \psi(v) \psi(v-u) \, dv.
\end{equation*}
A change of variables gives $\rho(-u) = \rho(u)$. Moreover,
\begin{equation*}
    \rho \in C_0^\infty(\R), \qquad \rho \geq 0, \qquad \supp \rho \subset [-1,1], \qquad \int_{\R} \rho(u) du = 1.
\end{equation*}
With the convention \eqref{fourier transform convention} of the Fourier transform, Fubini's theorem gives $\widehat \rho(\xi) = |\widehat \psi(\xi)|^2$. In particular,
\begin{equation}\label{rho positivity}
    0 \leq \widehat\rho(\xi) \leq 1, \quad \widehat \rho(0) = 1.
\end{equation}
Moreover, since $\rho \in C_0^\infty(\R)$, repeated integration by parts gives 
\begin{equation}\label{eq:rho-decay}
    |\widehat\rho(\xi)| \leq C_L(1 + |\xi|)^{-L}
\end{equation}
for every integer $L \geq 1$. For $R \geq 2$, put
\begin{equation*}
    \rho_R(u) \defeq R\rho(Ru), \quad b_{m,n}(R) \defeq \frac{\widehat\rho(2r_{m,n}/R)}{r_{m,n}^p},
\end{equation*}
and define
\begin{equation*}
    G_R(t) \defeq \sum_{m,n \geq 1}b_{m,n}(R) \cos(2tr_{m,n} - \theta_\gamma).
\end{equation*}
For each $R$, \eqref{eq:rho-decay} implies that $G_R(t)$ converges absolutely and uniformly for $t \in \R$. For $t > 1 + 1/R$, we claim that
\begin{equation}\label{shifted averaging}
    G_R(t) = \int_{\R}\rho_R(v)S_{\gamma,s}(t-v)\,dv.
\end{equation}
Indeed, since $\rho_R$ is even and real-valued, $\cos x = \operatorname{Re}(e^{ix})$ implies
\begin{equation*}
    \begin{split}
        \int_{\mathbb R}\rho_R(v)\cos(2(t-v)r-\theta_\gamma)\,dv &= \operatorname{Re}\bigg(\int_\R \rho_R(v) e^{i(2(t-v)r - \theta_\gamma)} dv\bigg) \\
        &= \operatorname{Re}(e^{i(2tr - \theta_\gamma)} \widehat{\rho_R}(2r)) = \widehat\rho(2r/R)\cos(2tr-\theta_\gamma),
    \end{split}
\end{equation*}
which gives 
\begin{equation*}
    \sum_{1 \leq m,n \leq N}b_{m,n}(R) \cos(2tr_{m,n}-\theta_\gamma) = \int_\R \rho_R(v) S_{\gamma,s,N}(t-v) dv.
\end{equation*}
%Indeed, evenness of $\rho_R$ implies
%\begin{equation*}
    %\int_{\mathbb R}\rho_R(v)
        %\cos(2(t-v)r-\theta_\gamma)\,dv
    %=
    %\widehat\rho(2r/R)\cos(2tr-\theta_\gamma).
%\end{equation*}
Letting $N \to \infty$, the left-hand side converges to $G_R(t)$, since the defining series converges absolutely. On the right-hand side, since $\supp \rho_R \subseteq [-1/R, 1/R]$, the arguments $t - v$ lie in the compact interval $[t - 1/R, t + 1/R] \subseteq [1,\infty)$, where the partial sums $S_{\gamma,s,N}$ converge uniformly by \eqref{eq: partial sums converge uniformly}. Since $\rho_R \in L^1(\R)$, \eqref{shifted averaging} follows from the dominated convergence theorem.
%On the right-hand side, the arguments $t-v$ lie in the compact interval $[t - 1/R, t + 1/R]$, where, as proved in \eqref{eq: partial sums converge uniformly}, the partial sums are uniformly bounded. The dominated convergence theorem then confirms \eqref{shifted averaging}. 

Next, set
\begin{equation*}
    A_R \defeq \sum_{m,n\geq1}b_{m,n}(R), \quad \alpha_\gamma \defeq -\cos\theta_\gamma.
\end{equation*}
Note that $G_R(0) = -\alpha_\gamma A_R$. Since $3\pi/4 < \theta_\gamma \leq \pi$, we have
$\alpha_\gamma > 0$.
%\begin{equation}\label{eq:GR-zero}
    %G_R(0)=-\alpha_\gamma A_R.
%\end{equation}
Define
\begin{equation*}
    h_\gamma(R) \defeq
    \begin{cases}
        R^{1/2-\gamma}, &0 < \gamma < \frac{1}{2}, \\
        \log R, &\gamma = \frac{1}{2}.
    \end{cases}
\end{equation*}
By \eqref{rho positivity} and continuity, there is an $\eta > 0$ such that
\begin{equation*}
    \widehat \rho(2r/R) \geq\frac{1}{2} \quad \text{when } r \leq \eta R.
\end{equation*}
Hence
\begin{equation*}
    A_R \geq \frac{1}{2}\sum_{r_{m,n} \leq \eta R}r_{m,n}^{-p}.
\end{equation*}
Furthermore, since $r_{m,n} \lesssim \max(m,n)$, there is a $c_0 > 0$ such that
\begin{equation*}
    \frac{1}{2}\sum_{r_{m,n} \leq \eta R}r_{m,n}^{-p} \gtrsim \sum_{\max(m,n) \leq c_0R} \max(m,n)^{-p}.
\end{equation*}
Each shell $\{(m, n) \in \NN^2: \max(m,n) = k\}$ contains exactly $2k - 1$ pairs. It follows that for sufficiently large $R$,
%Each shell $\{(m, n) \in \NN^2: \max(m,n)= k \}$, contains $2k - 1$ pairs. Furthermore, since $r_{m,n} \asymp \max(m,n)$, there is $c_0 > 0$ so that all the shells $k\leq c_0R$ lie in $r_{m,n}\leq\eta R$. It follows that for sufficiently large $R$,
%\begin{equation}\label{shifted mass lower}
    %A_R \geq c\sum_{k \leq c_0R} (2k-1)k^{-p} \geq c\,h_\gamma(R).
%\end{equation}
\begin{equation} \label{shifted mass lower}
    A_R \gtrsim \sum_{k \leq c_0R} (2k-1)k^{-p} \gtrsim \sum_{k \leq c_0R} k^{-1/2-\gamma} \gtrsim \,h_\gamma(R).
\end{equation}
For $K \geq 1$, define the tail
\begin{equation*}
    \mathcal T_R(K) \defeq \sum_{r_{m,n} > KR} b_{m,n}(R).
\end{equation*}
Taking $L = 4$ in \eqref{eq:rho-decay},
\begin{equation*}
    \mathcal T_R(K) \lesssim R^4\sum_{r_{m,n} > KR} r_{m,n}^{-p-4}.
\end{equation*}
To estimate this, partition pairs $(m,n)$ into
\begin{equation*}
    \mathcal A_j \defeq \{(m,n) \in \NN^2: 2^jKR < r_{m,n} \leq 2^{j+1}KR\}, \quad j \geq 0.
\end{equation*}
If $r_{m,n} \leq u$, then 
\begin{equation*}
    m = s + \frac{1}{a}\sqrt{a^2(m - s)^2} \leq s + \frac{u}{a}
\end{equation*}
%$m \leq s + u/a$ 
and similarly $n \leq s + u/b$. Thus, for $u \geq 1$,
\begin{equation} \label{eq:shifted-lattice-count}
    \#\{(m,n) \in \NN^2 : r_{m,n} \leq u\} \leq \bigg(s + \frac{u}{a}\bigg)\bigg(s + \frac{u}{b}\bigg) \lesssim u^2.
\end{equation}
In particular,
\begin{equation*}
    \#\mathcal A_j\lesssim (2^jKR)^2.
\end{equation*}
Thus,
\begin{equation} \label{shifted tail bound}
    \begin{split}
        \mathcal T_R(K) &\lesssim R^4 \sum_{j = 0}^\infty \sum_{(m,n) \in \mathcal{A}_j} r_{m,n}^{-p-4} \\
        &\lesssim R^4 \sum_{j=0}^\infty (2^jKR)^2(2^jKR)^{-p-4} \\
        &\lesssim R^{2-p}K^{-p-2} \lesssim h_\gamma(R)K^{-p-2}.
    \end{split}
\end{equation}
%Therefore
%\begin{equation*}
    %\sum_{r_{m,n}>KR}r_{m,n}^{-p-4} \lesssim \sum_{j=0}^\infty (2^jKR)^2(2^jKR)^{-p-4} \lesssim (KR)^{-p-2},
%\end{equation*}
%and consequently
%\begin{equation} \label{shifted tail bound}
    %\mathcal T_R(K) \lesssim R^{2-p}K^{-p-2}.
%\end{equation}
%For $0 < \gamma < \frac{1}{2}$, we have $R^{2 - p} = h_\gamma(R)$. At $\gamma = \frac{1}{2}$, we have $R^{2 - p} = 1 \lesssim \log R$. 
Combining \eqref{shifted mass lower} and \eqref{shifted tail bound}, we may fix a sufficiently large $K$ such that
\begin{equation}\label{eq:shifted-tail-small}
    \mathcal T_R(K) \leq \frac{\alpha_\gamma}{8}A_R
\end{equation}
for $R$ large enough.

Next, enumerate the radii satisfying $r_{m,n} \leq KR$ as $r_1, \dots, r_{D_R}$, allowing repetitions. By \eqref{eq:shifted-lattice-count},
\begin{equation} \label{shifted frequency count}
    D_R\lesssim K^2R^2.
\end{equation}
Fix a $Q \in \NN$ satisfying
\begin{equation*}
    \frac{2\pi}{Q} \leq \frac{\alpha_\gamma}{4}.
\end{equation*}
Applying Dirichlet's simultaneous approximation theorem \cite[Chapter II, Theorem 1A]{schmidt80} to the numbers
\begin{equation*}
    \beta_\ell \defeq \frac{Rr_\ell}{\pi}, \quad 1 \leq \ell \leq D_R,
\end{equation*}
gives an integer $q_R$ and integers $z_1, \dots, z_{D_R}$ such that
\begin{equation} \label{eq:shifted-dirichlet}
    1 \leq q_R \leq Q^{D_R}, \qquad |q_R\beta_\ell - z_\ell| \leq \frac{1}{Q}, \qquad 1 \leq \ell \leq D_R.
\end{equation}
Put $t_R \defeq Rq_R$. Then
\begin{equation} \label{eq:shifted-time-range}
    R \leq t_R \leq RQ^{D_R},
\end{equation}
and
\begin{equation} \label{eq: tR bound}
    |2t_Rr_\ell - 2\pi z_\ell| \leq \frac{2\pi}{Q}, \qquad 1 \leq \ell \leq D_R.
\end{equation}
By periodicity and the inequality $|\cos x - \cos y| \leq |x - y|$,
\begin{equation*}
    |\cos(2t_Rr_\ell - \theta_\gamma) - \cos\theta_\gamma| \leq \frac{2\pi}{Q}, \qquad 1 \leq \ell \leq D_R.
\end{equation*}
Adding and subtracting $G_R(0)$, and splitting the sum according to  $r_{m,n} \le KR$ or $r_{m,n} > KR$, we obtain
\begin{equation*}
    \begin{split}
        G_R(t_R) &= G_R(0) + \sum_{m,n \geq 1}b_{m,n}(R) \cos(2t_Rr_{m,n} - \theta_\gamma) - G_R(0) \\
        &= -\alpha_\gamma A_R \\
        &+ \sum_{r_{m,n}\leq KR}b_{m,n}(R) \bigl(\cos(2t_Rr_{m,n}-\theta_\gamma) -\cos\theta_\gamma\bigr) \\
        &+ \sum_{r_{m,n}>KR}b_{m,n}(R) \bigl(\cos(2t_Rr_{m,n}-\theta_\gamma)-\cos\theta_\gamma\bigr) \\
        &\leq -\alpha_\gamma A_R + \frac{2\pi}{Q}A_R + 2\mathcal T_R(K) \\
        &\leq \bigg(-\alpha_\gamma + \frac{\alpha_\gamma}{4} + \frac{\alpha_\gamma}{4}\bigg)A_R = -\frac{\alpha_\gamma}{2}A_R. 
    \end{split}
\end{equation*}
Since
\begin{equation*}
    t_R \geq R > 1 + \frac{1}{R},
\end{equation*}
we may apply \eqref{shifted averaging}. This gives
\begin{equation*}
    \int_{\R}\rho_R(v)S_{\gamma,s}(t_R - v)\,dv \leq -\frac{\alpha_\gamma}{2}A_R.
\end{equation*}
Since $\rho_R \geq 0$, $\int_\R \rho_R = 1$, and $\operatorname{supp}\rho_R \subseteq [-1/R,1/R]$, there is some $u_R \in [t_R - 1/R, t_R + 1/R]$ such that
\begin{equation*}
    S_{\gamma,s}(u_R) \leq -\frac{\alpha_\gamma}{4}A_R.
\end{equation*}
Combining this with \eqref{shifted mass lower}, we obtain
\begin{equation} \label{eq: S gamma s uR bound}
    S_{\gamma,s}(u_R) \leq -c\,h_\gamma(R),
\end{equation}
for some $c > 0$ independent of $R$. %Otherwise the average could not be at most $-\alpha_\gamma A_R/2$. 
Now let $E_R \defeq u_R^2$. From \eqref{eq:shifted-time-range} and $u_R \in [t_R - 1/R, t_R + 1/R]$,
\begin{equation*}
    \bigg(R - \frac{1}{R}\bigg)^2 \leq E_R \leq \bigg(RQ^{D_R} + \frac{1}{R}\bigg)^2.
\end{equation*}
Hence
\begin{equation*}
    \frac{R^2}{4} \leq E_R \leq 4R^2Q^{2D_R}.
\end{equation*}
%\begin{equation*}
    %R - \frac{1}{R} \leq u_R \leq RQ^{D_R} + \frac{1}{R}.
%\end{equation*}
%Thus, for sufficiently large $R$,
%\begin{equation*}
    %\frac{R^2}{4} \leq E_R \leq 4R^2Q^{2D_R}.
%\end{equation*}
In particular, $E_R \to \infty$, and \eqref{shifted frequency count} implies
\begin{equation*}
    \log E_R \lesssim R^2.
    %\log E_R
    %\leq\log4+2\log R+2D_R\log Q
    %\lesssim R^2.
\end{equation*}
It follows that there is a $c > 0$ such that
\begin{equation*}
    h_\gamma(R) \geq c
    \begin{cases}
        (\log E_R)^{1/4-\gamma/2}, & 0 < \gamma < \frac{1}{2}, \\
        \log\log E_R, & \gamma = \frac{1}{2},
    \end{cases}
\end{equation*}
for $R$ large enough. Combining this with \eqref{eq: S gamma s uR bound} and extracting an increasing sequence $E_j \to \infty$ completes the proof of \eqref{shifted double sum lower}.
%For $R$ large enough and $0 < \gamma < \frac{1}{2}$, we thus have
%\begin{equation*}
    %R^{1/2 - \gamma} \geq c(\log E_R)^{1/4-\gamma/2},
%\end{equation*}
%while for $\gamma = \frac{1}{2}$,
%\begin{equation*}
    %\log R \geq c \log\log E_R.
%\end{equation*}
%This completes the proof of \eqref{shifted double sum lower}.\\
\end{proof}
%\paragraph{The two double sums for the triangles.}
For the equilateral and hemi-equilateral triangles, the leading double-series contributions in \eqref{eq: Riesz mean equilateral asymptotic} and \eqref{eq: Riesz mean hemiequilateral asymptotic} both involve a linear combination of the double sums \eqref{eq: S gamma s}. Indeed, taking $a = \sqrt{3}$, $b = 1$, and $t = \ell\sqrt{3E}/2$, these contributions are given by
\begin{equation*}
    W_\gamma(t) \defeq 2S_{\gamma,0}(t)+2S_{\gamma,1/2}(t).
\end{equation*}
The upper bound \eqref{shifted double sum upper} follows immediately for $W_\gamma(t)$ by applying the corresponding bound to each sum. To establish unboundedness below, however, we must account for possible cancellation between the two double sums. For $s \in \{0,\frac{1}{2}\}$, let $G_{R,s}$, $A_{R,s}$, and $\mathcal{T}_{R,s}(K)$ denote the quantities $G_R$, $A_R$, and $\mathcal{T}_R(K)$, respectively, defined in the proof of Lemma \ref{shifted double sum bounds proposition}. Fix a single sufficiently large $K$ such that \eqref{eq:shifted-tail-small} holds for both $s = 0$ and $s = \frac{1}{2}$. Then
\begin{equation*}
    G_R^{W}(t) \defeq 2G_{R,0}(t) + 2G_{R,1/2}(t), \quad A_R^{W} \defeq 2A_{R,0} + 2A_{R,1/2}
\end{equation*}
satisfy
\begin{equation*}
    G_R^{W}(0) = -\alpha_\gamma A_R^{W}, \quad A_R^{W} \gtrsim h_\gamma(R),
\end{equation*}
where $\alpha_\gamma, h_\gamma$ are defined in the proof of Lemma \ref{shifted double sum bounds proposition}. Next, form a list consisting of all radii $r_{m,n} \leq KR$ corresponding to both $s = 0$ and $s = \frac{1}{2}$, with repetitions allowed. By \eqref{shifted frequency count}, the length of this combined list is still $O(R^2)$. Applying Dirichlet's simultaneous approximation theorem \cite[Chapter II, Theorem 1A]{schmidt80} to the numbers $Rr/\pi$ corresponding to these radii, we obtain a single $t_R = Rq_R$ satisfying estimates analogous to \eqref{eq: tR bound} for all radii in the combined list. Then
\begin{equation*}
    G_R^{W}(t_R) \leq \bigg(-\alpha_\gamma  + \frac{\alpha_\gamma}{4} + \frac{\alpha_\gamma}{4}\bigg)A_R^{W} = -\frac{\alpha_\gamma}{2}A_R^{W}.
\end{equation*}
Moreover,
\begin{equation*}
    G_R^{W}(t_R) = \int_{\R} \rho_R(v)W_\gamma(t_R - v)\,dv.
\end{equation*}
The remainder of the proof is unchanged. In particular, there is a $c > 0$ and an increasing sequence $E_j \to \infty$ such that $W_\gamma(\sqrt{E_j})$ satisfies the lower bound \eqref{shifted double sum lower}.

In the Riesz mean expansions \eqref{eq: Riesz mean equilateral asymptotic}-\eqref{eq: Riesz mean isosceles asymptotic}, the arguments of the double sums are fixed positive multiples of $\sqrt E$, which do not affect the logarithmic lower bounds, up to constants. Moreover, the additional single series in the leading coefficient functions converge absolutely and are uniformly bounded. Together with Lemma \ref{shifted double sum bounds proposition}, this shows that the leading coefficient functions $F_{\Omega,1}^\gamma$ are unbounded for all integrable polygons when $0 < \gamma \leq 1/2$.

\end{document}